\documentclass[thmsa,final,french,11pt]{article}
\usepackage{amsfonts}
\usepackage{amsmath}
\usepackage{amssymb}

\newtheorem{Th}{Theorem}
\newtheorem{Prop}[Th]{Proposition}

\newtheorem{Rk}{Remark}

\renewcommand{\phi}{\varphi}

\def\WW_#1{\boldsymbol{W}\!_#1}
\input{tcilatex}

\input tcilatex
\begin{document}

\title{Power variations for a class of Brown-Resnick processes }
\author{Christian Y. ROBERT\thanks{%
Universit\'{e} de Lyon, Universit\'{e} Lyon 1, Institut de Science Financi%
\`{e}re et d'Assurances, 50 Avenue Tony Garnier, F-69007 Lyon, France}}
\maketitle

\begin{abstract}
We consider the class of simple Brown-Resnick max-stable processes whose
spectral processes are continuous exponential martingales. We develop the
asymptotic theory for the realized power variations of these max-stable
processes, that is, sums of powers of absolute increments. We consider an
infill asymptotic setting, where the sampling frequency converges to zero
while the time span remains fixed. More specifically we obtain a biased
central limit theorem whose bias depend on the local times of the
differences between the logarithms of the underlying spectral processes. 

\strut

\end{abstract}

\section{Introduction}

In the two last decades there has been an increasing interest in limit
theory for power variations of stochastic processes because such functionals
are very important in analyzing the fine properties of the underlying model
and in statistical inference. Asymptotic theory for power variations of
various classes of stochastic processes has been intensively investigated in
the literature. We refer e.g. to Jacod and Protter (2011) for limit theory
for power variations of It\^{o} semimartingales, to Barndorff-Nielsen,
Corcuera and Podolskij (2009) for asymptotic results in the framework of
fractional Brownian motion and related processes.

In this paper we study the power variations for a class of max-stable
stochastic processes: the \textit{simple} Brown-Resnick max-stable processes
whose spectral processes are continuous exponential martingales. The
original process of this class was introduced in the seminal paper of Brown
and Resnick (1977) and it has been generalized to form a flexible family of
stationary max-stable processes based on Gaussian random fields by
Kabluchko, Schlather and de Haan (2009). The characterization of general
max-stable stochastic processes in $C[0,1]$, the space of continuous
functions on $[0,1]$ has been provided by Gin\'{e}, Hahn and Vatan (1990) at
the beginning of the nineties, while a decade later, de Haan and Lin (2001)
investigated the domain of attraction conditions.

To the best of our knowledge, our paper studies for the first time power
variations of a max-stable process. We provide several central limit
theorems in an infill asymptotic setting, i.e. where the sampling frequency
converges to zero while the time span remains fixed. It should however be
underlined that such an asymptotic approach has been considered for the
estimation of the integrated variance of a white noise process with a
positive and constant extreme value index by Einmahl, de Haan and Zhou
(2016). However such a process is not realistic for a large number of
applications for which the assumption of independent observations at any
high frequencies may appear as too strong. With their assumption they can
work with the values of the process (this is more or less equivalent to
consider triangular arrays of independent random variables), while we have
to consider the increments of the processes to mitigate the (strong) local
dependence between the observations.

The paper is organized as follows. Section 2 is devoted to presenting the
setting and to providing definitions and assumptions. Section 3 discusses
asymptotic results on the normalized power variations of the maximum of two
independent Brownian motions, of the original Brown-Resnick process, and of
a general max-stable process in our class of Brown-Resnick processes. All
proofs are gathered in Appendix.

\section{Setting, definitions and assumptions}

We consider a filtered probability space $(\Omega ,\mathcal{F},\left( 
\mathcal{F}_{t}\right) _{t\in \lbrack 0,1]},\mathbb{P})$. We denote by $%
C^{+}[0,1]:=\{f\in C[0,1]:f>0\}$ the space of continuous and positive
functions on $[0,1]$. Note that we equip $C[0,1]$ and $C^{+}[0,1]$ with the
supremum norm $\left\vert f\right\vert _{\infty }=\sup_{s\in \lbrack
0,1]}|f_{s}|$.

A stochastic process $\eta $ on $C^{+}[0,1]$ with nondegenerate marginals is
called \textit{simple }max-stable if for all positive integers $k$%
\begin{equation*}
\frac{1}{k}\bigvee_{i=1}^{k}\eta _{i}\overset{\mathcal{L}}{=}\eta
\end{equation*}%
where $\eta _{1},\eta _{2},...$ are independent and identically distributed
(iid) copies of the process $\eta $, $\mathbb{P}(\eta _{t}<1)=e^{-1}$ for
all $t\in \lbrack 0,1]$ (i.e. it has a standard Frechet distribution) and $%
\overset{\mathcal{L}}{=}$ means equality in distribution. By Corollary 9.4.5
in de Haan and Ferreira (2006), all simple max-stable processes in $%
C^{+}[0,1]$ can be generated in the following way. Consider a Poisson point
process, $(R_{i})_{i\geq 1}$, on $(0,\infty )$ with mean measure $dr/r^{2}$.
Further consider iid stochastic processes $V,V_{1},V_{2},...$ in $C^{+}[0,1]$
with $\mathbb{E[}V_{t}]=1$ for all $t\in \lbrack 0,1]$ and $\mathbb{E[}%
\sup_{t\in \lbrack 0,1]}V_{t}]<\infty $. Let the point process and the
sequence $V_{1},V_{2},...$ be independent. Then%
\begin{equation}
\eta =\bigvee_{i=1}^{\infty }R_{i}V_{i}  \label{SimpleMaxStable}
\end{equation}%
is a simple max-stable process. Conversely, each simple max-stable process
has such a representation (which is not unique). The process $V$ is called a 
\textit{spectral process} associated to $\eta $.

We now consider a sequence of iid random processes $\xi ,\xi _{1},\xi
_{2},...$ in $C[0,1]$. This sequence is said to belong to the domain of
attraction of the simple max-stable process $\eta $ if there exists a
sequence of non-random positive normalizing functions $\left( c_{n,t}\right)
_{t\in \lbrack 0,1]}$, $n=1,...$, such that%
\begin{equation*}
\left\{ c_{n,t}^{-1}\bigvee_{i=1}^{n}\xi _{i,t}\right\} _{t\in \lbrack 0,1]}%
\overset{\mathcal{L}}{\implies }\eta ,\qquad n\rightarrow \infty ,
\end{equation*}%
where $\overset{\mathcal{L}}{\implies }$ denotes the convergence in law in $%
C[0,1]$, see e.g. Theorem 9.2.1 in de Haan and Ferreira (2006). Let $%
S_{c}:=\{f\in C^{+}[0,1];\left\vert f\right\vert _{\infty }\geq c\}$ where $%
c>0$, and define a sequence of measures $\nu _{n}$, $n=1,2,...$, by%
\begin{equation*}
\nu _{n}(\cdot ):=n\Pr \left( n^{-1}\zeta \in \left( \cdot \right) \right)
\end{equation*}%
to $S_{c}$, where $\zeta _{t}:=(1-F_{t}\left( \xi _{t}\right) )^{-1}$ and $%
F_{t}(x)=\mathbb{P}\left( \xi _{t}\leq x\right) $. The sequence of measures $%
\nu _{n}$ weakly converges, as $n\rightarrow \infty $, to the restriction of
the so-called \textit{exponent measure} $\nu $ of $\eta $ to $S_{c}$ for
each $c>0$. It can be shown that the exponent measure coincides with the
distribution of the spectral process, i.e. $\nu (\cdot )=\mathbb{P}\left(
V\in \left( \cdot \right) \right) $ and that 
\begin{equation*}
\lim_{n\rightarrow \infty }\frac{c_{\left\lfloor nx\right\rfloor ,t}}{c_{n,t}%
}=x
\end{equation*}%
locally uniformly for $x\in (0,\infty )$ and uniformly for $t\in \lbrack
0,1] $.

In this paper, we will assume that $V$ is a continuous exponential
martingale defined by 
\begin{equation}
V_{t}=\exp \left\{ \int_{0}^{t}H_{s}\text{d}W_{s}-\frac{1}{2}%
\int_{0}^{t}H_{s}^{2}\text{d}s\right\} ,\quad t\in \lbrack 0,1],
\label{ExpMart}
\end{equation}%
where $H$ is a non-random H\"{o}lder function in $C^{+}[0,1]$ with exponent $%
\alpha >1/2$ and satisfying $\int_{0}^{1}H_{s}^{2}ds<\infty $ and $%
\inf_{s\in \lbrack 0,1]}H_{s}>c$ for some positive constant $c$, $%
\{W_{t},t\in \lbrack 0,1]\}$ is a $\mathcal{F}$-adapted standard Brownian
motion on $[0,1]$. All the processes $V_{i}$, $i\geq 1$, in Eq. $\left( \ref%
{SimpleMaxStable}\right) $ are $\mathcal{F}$-adapted. We call the family of
processes $\eta $ associated with $V$ the family of Brown-Resnick processes
because when $H_{t}=\sigma >0$ for all $t\in \lbrack 0,1]$, $\eta $ is the
stationary max-stable process introduced in Brown and Resnick (1977).

\section{Asymptotic behaviors of normalized power variations of
Brown-Resnick processes}

The increments of a stochastic process $X$ over the equi-spaced grid with
mesh $1/n$ of $[0,1]$ are denoted by%
\begin{equation*}
\Delta _{i}^{n}X=X_{i/n}-X_{(i-1)/n}\quad i=1,..,n.
\end{equation*}%
The normalized power variation of order $p\geq 1$ of $X$ is defined by 
\begin{equation*}
B\left( p,X\right) _{t}^{n}=n^{p/2-1}\sum_{i=1}^{\left\lfloor
nt\right\rfloor -1}|\Delta _{i}^{n}X|^{p}.
\end{equation*}%
In this section we discuss the asymptotic behavior of $B\left( p,X\right)
_{t}^{n}$ for several stochastic processes. We begin with the maximum of two
independent Brownian motions, then with the logarithm of the original Brown
Resnick process for which $H_{t}=\sigma >0$ for $t\in \lbrack 0,1]$, and
finally with the logarithm of $\eta $.

We denote by $\overset{u.c.p}{\Longrightarrow }$ the convergence in
probability, uniform over each compact interval in $[0,1]$. We also need to
recall the notion of stable convergence in law, which was introduced in R%
\'{e}nyi (1963). Let $Z_{n}$ be a sequence of $E$-valued random variables
defined on the same probability space $(\Omega ,\mathcal{F},\mathbb{P})$.
Let $Z$ be an $E$-valued random variable defined on an extension, $(\tilde{%
\Omega},\mathcal{\tilde{F}},\mathbb{\tilde{P}})$. We then say that $Z_{n}$
converges $\mathcal{F}$-stably to $Z$ (and write $Z_{n}\overset{\mathcal{L}-s%
}{\implies }Z$) if%
\begin{equation*}
\lim_{n\rightarrow \infty }\mathbb{E}[Uf(Z_{n})]=\mathbb{\tilde{E}}[Uf(Z)]
\end{equation*}%
for all bounded continuous functions $f$ on $E$ and all bounded $\mathcal{F}$%
-measurable random variables $U$. This notion of convergence is stronger
than convergence in law, but weaker than convergence in probability. We
refer to Jacod and Protter (2011) for a detailed exposition of this last
type of convergence.

\subsection{Normalized power variations of two independent Brownian motions}

We here consider the case of the maximum of two Brownian motions $W_{1}\vee
W_{2}=\{W_{1,t}\vee W_{2,t},t\in \lbrack 0,1]\}$, where $W_{1}=\{W_{1,t},t%
\in \lbrack 0,1]\}$ and $W_{2}=\{W_{2,t},t\in \lbrack 0,1]\}$ are two
independent Brownian motions defined on $(\Omega ,\mathcal{F},\left( 
\mathcal{F}_{t}\right) _{t\in \lbrack 0,1]},\mathbb{P})$. Let us recall that
(see e.g. p. 10 in Jacod and Protter (2011)) that%
\begin{equation*}
B\left( p,W_{1}\right) _{t}^{n}\overset{u.c.p}{\Longrightarrow }m_{p}t\text{%
\qquad and\qquad }\sqrt{n}\left( B\left( p,W_{1}\right)
_{t}^{n}-m_{p}t\right) \overset{\mathcal{L}}{\implies }\tilde{X}_{t}
\end{equation*}%
where $m_{p}$ is the expectation of the $p$-th moment of the absolute value
of a standard Gaussian random variable, and $\tilde{X}$ is a continuous
centered Gaussian martingale with variance $(m_{2p}-m_{p}^{2})t$.

Let us denote by $(x)_{+}$ the positive part of a real $x$ and let $%
W_{2\backslash 1}=W_{2}-W_{1}$. Since $W_{1}\vee W_{2}=W_{1}+(W_{2\backslash
1})_{+}$, we deduce by Tanaka's formula that%
\begin{equation*}
W_{1,t}\vee W_{2,t}=W_{1,t}+\int_{0}^{t}\mathbb{I}_{\{W_{2\backslash
1,t}>0\}}dW_{2\backslash 1,t}+\frac{1}{2}L_{_{W_{2\backslash
1}},t}^{0},\quad t\in \lbrack 0,1],
\end{equation*}%
where $L_{_{W_{2\backslash 1}},t}^{0}$ is the local time of $W_{2\backslash
1}$ at time $t$ and level $0$. As a consequence $W_{1}\vee W_{2}$ is not an
Ito semi-martingale (since its predictable part of finite variation is not
absolutely continuous with respect to the Lebesgue measure). The asymptotic
results of functionals of normalized increments of a semi-martingale are
often obtained under the assumption that the semi-martingale is an Ito
semi-martingale (see e.g. Section 3.4.2 and 5.3 in Jacod and Protter
(2011)). Therefore the results given in Jacod and Protter (2011) can not be
used directly in our case.

Let $f$ be a real measurable function. By partitioning on the positive and
negative values of $W_{2\backslash 1,(i-1)/n}$ and $W_{2\backslash 1,i/n}$,
we have%
\begin{eqnarray*}
&&f\left( \sqrt{n}\Delta _{i}^{n}(W_{1}\vee W_{2})\right) \\
&=&f(\sqrt{n}\Delta _{i}^{n}W_{1})\mathbb{I}_{\{W_{2\backslash
1,(i-1)/n}<0,W_{2\backslash 1,i/n}<0\}}+f(\sqrt{n}\Delta _{i}^{n}W_{2})%
\mathbb{I}_{\{W_{2\backslash 1,(i-1)/n}>0,W_{2\backslash 1,i/n}>0\}} \\
&&+f\left( \sqrt{n}\left( W_{2,i/n}-W_{1,(i-1)/n}\right) \right) \mathbb{I}%
_{\{W_{2\backslash 1,(i-1)/n}\leq 0,W_{2\backslash 1,i/n}\geq 0\}} \\
&&+f\left( \sqrt{n}\left( W_{1,i/n}-W_{2,(i-1)/n}\right) \right) \mathbb{I}%
_{\{W_{2\backslash 1,(i-1)/n}\geq 0,W_{2\backslash 1,i/n}\leq 0\}}.
\end{eqnarray*}%
It follows that%
\begin{eqnarray}
f\left( \sqrt{n}\Delta _{i}^{n}(W_{1}\vee W_{2})\right) &=&f\left( \sqrt{n}%
\Delta _{i}^{n}W_{1}\mathbb{I}_{\{W_{2\backslash 1,(i-1)/n}<0\}}+\sqrt{n}%
\Delta _{i}^{n}W_{2}\mathbb{I}_{\{W_{2\backslash 1,(i-1)/n}>0\}}\right) 
\notag \\
&&+\Psi _{f}(\sqrt{n}\Delta _{i}^{n}W_{1},\sqrt{n}\Delta _{i}^{n}W_{2},\sqrt{%
n}W_{2\backslash 1,(i-1)/n})  \label{f_MaxB}
\end{eqnarray}%
where 
\begin{equation*}
\Psi _{f}\left( x,y,w\right) =(f(y+w)-f(x))\mathbb{I}_{\{x-y\leq w\leq
0\}}+(f(x-w)-f(y))\mathbb{I}_{\{0\leq w\leq x-y\}}.
\end{equation*}%
One can remark that 
\begin{equation*}
\Delta _{i}^{n}W_{1}\mathbb{I}_{\{W_{2\backslash 1,(i-1)/n}<0\}}+\Delta
_{i}^{n}W_{2}\mathbb{I}_{\{W_{2\backslash 1,(i-1)/n}>0\}}
\end{equation*}%
has the same distribution as $\Delta _{i}^{n}W_{1}$ or $\Delta _{i}^{n}W_{2}$
and is independent of $\sigma \left( W_{1,t},0\leq t\leq (i-1)/n\right) 
\mathcal{\ }$and of $\sigma \left( W_{2,t},0\leq t\leq (i-1)/n\right) $.

Let us now define%
\begin{equation*}
\varphi _{p}\left( w\right) =\int_{\mathbb{R}^{2}}\Psi _{|\cdot |^{p}}(x,y,w)%
\frac{1}{2\pi }e^{-(x^{2}+y^{2})/2}\text{d}x\text{d}y.
\end{equation*}%
We can remark that 
\begin{equation*}
\varphi _{p}\left( \sqrt{n}W_{2\backslash 1,(i-1)/n}\right) =\mathbb{E}\left[
\left. \Psi _{|\cdot |^{p}}(\sqrt{n}\Delta _{i}^{n}W_{1},\sqrt{n}\Delta
_{i}^{n}W_{2},\sqrt{n}W_{2\backslash 1,(i-1)/n})\right\vert \mathcal{F}%
_{(i-1)/n}\right] .
\end{equation*}%
Since $\int |\varphi _{p}\left( w\right) |$d$w<\infty $ for any $p\geq 1$,
we can deduce from Theorem 1.1 in Jacod (1998), that%
\begin{equation*}
\frac{1}{\sqrt{n}}\sum_{i=1}^{\left\lfloor nt\right\rfloor -1}\varphi
_{p}\left( \sqrt{n}W_{2\backslash 1,(i-1)/n}\right) \overset{u.c.p}{%
\Longrightarrow }\frac{1}{2}\lambda (\varphi _{p})L_{W_{2\backslash
1},t}^{0}.
\end{equation*}%
where $\lambda (\varphi _{p})=\int \varphi _{p}\left( w\right) $d$w$.

We now state our first result.

\begin{Prop}
\label{TCLMaxB}As $n\rightarrow \infty $, we have%
\begin{equation*}
B\left( p,W_{1}\vee W_{2}\right) _{t}^{n}\overset{u.c.p}{\Longrightarrow }%
m_{p}t\text{\qquad and\qquad }\sqrt{n}\left( B\left( p,W_{1}\vee
W_{2}\right) _{t}^{n}-m_{p}t\right) \overset{\mathcal{L}-s}{\implies }\tilde{%
X}_{1,t},
\end{equation*}%
where $\tilde{X}_{1}$ is a process defined on an extension $(\Omega ,%
\mathcal{\tilde{F}},(\mathcal{\tilde{F})}_{t\geq 0},\mathbb{P})$ of $(\Omega
,\mathcal{F},\left( \mathcal{F}\right) _{t\geq 0},\mathbb{P})$, which
conditionally on $\mathcal{F}$ is a continuous Gaussian process, with
independent increments, and whose mean and variance are given respectively by%
\begin{equation*}
\mathbb{\tilde{E}}[\tilde{X}_{1,t}|\mathcal{F]}=\frac{1}{2}\lambda (\varphi
_{p})L_{W_{2\backslash 1},t}^{0}\quad \text{and}\quad \mathbb{\tilde{V}}[%
\tilde{X}_{1,t}|\mathcal{F]}=(m_{2p}-m_{p}^{2})t.
\end{equation*}
\end{Prop}

We observe that, contrary to the Brownian case, the asymptotic convergence
of 
\begin{equation*}
\sqrt{n}\left( B\left( p,W_{1}\vee W_{2}\right) _{t}^{n}-m_{p}t\right)
\end{equation*}%
needs the stable convergence in law because of the additional term $\lambda
(\varphi _{p})L_{W_{2\backslash 1},t}^{0}/2$ as the (conditional) mean of $%
\tilde{X}_{1,t}$.

\subsection{Normalized power variations of the logarithm of the Brown
Resnick processes}

We first consider the case of the original Brown Resnick process for which $%
H_{t}=\sigma >0$ for $t\in \lbrack 0,1]$. We have 
\begin{equation*}
\log \eta _{t}=\bigvee_{i=1}^{\infty }\left( \log R_{i}+\sigma W_{i,t}-\frac{%
1}{2}\sigma ^{2}t\right) =\bigvee_{i=1}^{\infty }\left( \log R_{i}+\sigma
W_{i,t}\right) -\frac{1}{2}\sigma ^{2}t
\end{equation*}%
where $(R_{i})_{i\geq 1}$ is a $\mathcal{F}_{0}$-Poisson point process, with
mean measure $dr/r^{2}$, and $W_{1},W_{2},...$ are independent $\left( 
\mathcal{F}_{t}\right) $-Brownian motions. It is well known (see e.g.
Kabluchko, Schlather and de Haan (2009)) that $\log \eta $ is a stationary
process.

Let us study the distribution of its normalized increments 
\begin{equation*}
U_{i}^{n}=\sqrt{n}\sigma ^{-1}\Delta _{i}^{n}\log \eta ,\quad i=1,...,n.
\end{equation*}%
The following proposition provides the conditional and marginal
distributions of these increments and allows to deduce that they have
asymptotically a standard Gaussian distribution.

\begin{Prop}
\label{PropIncMAS}Let $u\in \mathbb{R}$. The conditional distribution of $%
U_{i}^{n}$ given $\eta _{(i-1)/n}=\eta $ is characterized by%
\begin{eqnarray*}
&&\Pr \left( \left. U_{i}^{n}\leq u\right\vert \eta _{(i-1)/n}=\eta \right)
\\
&=&\exp \left( -\frac{1}{\eta }\left[ e^{-\sigma u/\sqrt{n}}\Phi \left( -u+%
\frac{\sigma }{2\sqrt{n}}\right) -\Phi \left( -u-\frac{\sigma }{2\sqrt{n}}%
\right) \right] \right) \Phi \left( u+\frac{\sigma }{2\sqrt{n}}\right) ,
\end{eqnarray*}%
and its marginal distribution by%
\begin{equation*}
\Pr \left( U_{i}^{n}\leq u\right) =\frac{\Phi \left( u+\sigma /(2\sqrt{n}%
)\right) }{\Phi \left( u+\sigma /(2\sqrt{n})\right) +e^{-\sigma u/\sqrt{n}%
}\left( 1-\Phi \left( u-\sigma /(2\sqrt{n})\right) \right) },
\end{equation*}%
where $\Phi $ is the cumulative distribution function of the standard
Gaussian distribution. Moreover, we have, for any $p\geq 1$, 
\begin{equation*}
\lim_{n\rightarrow \infty }\sqrt{n}\left( \mathbb{E}\left[ |U_{i}^{n}|^{p}%
\right] -m_{p}\right) =2p\int_{0}^{\infty }u^{p-1}\varphi (u)\left[ 1/2-\bar{%
\Phi}\left( u\right) -u\bar{\Phi}\left( u\right) \Phi \left( u\right)
/\varphi (u)\right] du,
\end{equation*}%
where $\varphi $ and $\bar{\Phi}$ are respectively the probability density
function and the survival distribution function of the standard Gaussian
distribution.
\end{Prop}

One can observe that%
\begin{equation*}
\Pr \left( U_{i}^{n}\leq u\right) +\Pr \left( U_{i}^{n}\leq -u\right)
=1,\quad u\in \mathbb{R},
\end{equation*}%
and we therefore conclude that $U_{i}^{n}$ has a symmetric distribution.
Moreover it is easily derived that the distribution of $U_{i}^{n}$ converges
to a standard Gaussian distribution. The rate of convergence of $\mathbb{E}%
\left[ |U_{i}^{n}|^{p}\right] $ to $m_{p}$ is however relatively slow ($1/%
\sqrt{n}$) and this has for consequence that 
\begin{equation*}
\lim_{n\rightarrow \infty }\sqrt{n}\left( \mathbb{E}\left[ B\left( p,\log
\eta \right) _{t}^{n}\right] -m_{p}\sigma ^{p}t\right) =2p\sigma
^{p}\int_{0}^{\infty }u^{p-1}\varphi (u)\left[ 1/2-\bar{\Phi}\left( u\right)
-u\bar{\Phi}\left( u\right) \Phi \left( u\right) /\varphi (u)\right] du.
\end{equation*}%
Therefore an asymptotic bias is expected in the limit of $\sqrt{n}\left(
B\left( p,\log \eta \right) _{t}^{n}-m_{p}\sigma ^{p}t\right) .$

Let us now introduce some notation. For $i,j,k\geq 1$, let%
\begin{equation*}
Z_{i,t}=\log R_{i}+\sigma W_{i,t}\quad \text{and}\quad Z_{k\backslash
j,t}=Z_{k,t}-Z_{j,t}.
\end{equation*}%
Let $f$ be a real measurable function. By partitioning on the values of $j$
and $k$ for which $\vee _{l\geq 1}Z_{l,(i-1)/n}=Z_{j,(i-1)/n}$ and $\vee
_{l\geq 1}Z_{l,i/n}=Z_{k,i/n}$, we have%
\begin{eqnarray*}
&&f\left( \sqrt{n}\Delta _{i}^{n}\log \eta +\frac{1}{2\sqrt{n}}\sigma
^{2}\right) \\
&=&\sum_{j\geq 1}f\left( \sqrt{n}\sigma \Delta _{i}^{n}W_{j}\right) \mathbb{I%
}_{\{\vee _{l\geq 1}Z_{l,(i-1)/n}=Z_{j,(i-1)/n},\vee _{l\geq
1}Z_{l,i/n}=Z_{j,i/n}\}} \\
&&+\sum_{j\geq 1}\sum_{k\neq j}f\left( \sqrt{n}Z_{k\backslash j,i/n}\right) 
\mathbb{I}_{\{\vee _{l\geq 1}Z_{l,(i-1)/n}=Z_{j,(i-1)/n},\vee _{l\geq
1}Z_{l,i/n}=Z_{k,i/n}\}.}
\end{eqnarray*}%
Let $\Psi _{f,\sigma }\left( x,y,w\right) =\Psi _{f}\left( \sigma x,\sigma
y,w\right) $. We have%
\begin{eqnarray*}
&&f\left( \sqrt{n}\Delta _{i}^{n}\log \eta +\frac{1}{2\sqrt{n}}\sigma
^{2}\right) \\
&=&f\left( \sigma \sum_{j\geq 1}\sqrt{n}\Delta _{i}^{n}W_{j}\mathbb{I}%
_{\{\vee _{l\geq 1}Z_{l,(i-1)/n}=Z_{j,(i-1)/n}\}}\right) \\
&&+\sum_{j\geq 1}\sum_{k>j}\Psi _{f,\sigma }(\sqrt{n}\Delta _{i}^{n}W_{j},%
\sqrt{n}\Delta _{i}^{n}W_{k},\sqrt{n}Z_{k\backslash j,(i-1)/n})\mathbb{I}%
_{\{\vee _{l\neq j,k}Z_{l,(i-1)/n}\leq \wedge _{l=j,k}Z_{l,(i-1)/n}\}} \\
&&+H_{f}^{n},
\end{eqnarray*}%
where 
\begin{eqnarray*}
H_{f}^{n} &=&\sum_{j\geq 1}\sum_{k\neq j}\Psi _{f,\sigma }^{<}(\sqrt{n}%
\Delta _{i}^{n}W_{j},\sqrt{n}\Delta _{i}^{n}W_{k},\sqrt{n}Z_{k\backslash
j,(i-1)/n})\times \\
&&\dprod\limits_{l\neq j,k}\mathbb{I}_{\{Z_{l\backslash j,(i-1)/n}\leq
0\}}\times \left( \dprod\limits_{l\neq k,j}\mathbb{I}_{\{Z_{l\backslash
k,i/n}\leq 0\}}-\dprod\limits_{l\neq k,j}\mathbb{I}_{\{Z_{l\backslash
j,(i-1)/n}\leq 0\}}\right) ,
\end{eqnarray*}%
with%
\begin{equation*}
\Psi _{f,\sigma }^{<}\left( x,y,w\right) =(f(\sigma y+w)-f(\sigma x))\mathbb{%
I}_{\{\sigma \left( x-y\right) \leq w\leq 0\}}.
\end{equation*}%
One can observe that 
\begin{equation*}
\sum_{j\geq 1}\Delta _{i}^{n}W_{j}\mathbb{I}_{\{\vee _{l\geq
1}Z_{l,(i-1)/n}=Z_{j,(i-1)/n}\}}
\end{equation*}%
has the same distribution as $\Delta _{i}^{n}W_{j}$, $j\geq 1$, and are
independent of the $\sigma \left( W_{j,t},0\leq t\leq (i-1)/n\right) $, $%
j\geq 1$.

Let us define%
\begin{equation*}
\varphi _{p,\sigma }\left( w\right) =\int_{\mathbb{R}^{2}}\Psi _{|\cdot
|^{p},\sigma }(x,y,w)\frac{1}{2\pi }e^{-(x^{2}+y^{2})/2}\text{d}x\text{d}y.
\end{equation*}%
We can deduce from a simple modification of Theorem 1.1 in Jacod (1998),
that, for $j\geq 1$ and $k>j$,%
\begin{eqnarray*}
&&\frac{1}{\sqrt{n}}\sum_{i=1}^{\left\lfloor nt\right\rfloor -1}\mathbb{E}%
\left[ \left. \Psi _{|\cdot |^{p},\sigma }(\sqrt{n}\Delta _{i}^{n}W_{j},%
\sqrt{n}\Delta _{i}^{n}W_{k},\sqrt{n}Z_{k\backslash j,(i-1)/n})\right\vert 
\mathcal{F}_{(i-1)/n}\right] \mathbb{I}_{\{\vee _{l\neq
j,k}Z_{l,(i-1)/n}\leq \wedge _{l=j,k}Z_{l,(i-1)/n}\}} \\
&&\overset{u.c.p}{\Longrightarrow }\frac{1}{2\sigma ^{2}}\lambda (\varphi
_{p,\sigma })\int_{0}^{t}\mathbb{I}_{\{\wedge _{l=j,k}Z_{l,s}>\vee _{l\neq
j,k}Z_{l,s}\}}dL_{Z_{k\backslash j},s}^{0},
\end{eqnarray*}%
where $L_{Z_{k\backslash j},s}^{0}$ is the local time of $Z_{k\backslash j}$
at time $t$ and level $0$. We can now state our result on the convergence of 
$B\left( p,\log \eta \right) _{t}^{n}$.

\begin{Prop}
\label{TCLSimpMS} Assume that $H_{t}=\sigma >0$ for $t\in \lbrack 0,1]$. As $%
n\rightarrow \infty $, we have for any integer $p\geq 1$%
\begin{equation*}
B\left( p,\log \eta \right) _{t}^{n}\overset{u.c.p}{\Longrightarrow }%
m_{p}\sigma ^{p}t\quad \text{and}\quad \sqrt{n}\left( B\left( p,\log \eta
\right) _{t}^{n}-m_{p}\sigma ^{p}t\right) \overset{\mathcal{L}-s}{\implies }%
\tilde{X}_{2,t},
\end{equation*}%
where $\tilde{X}_{2}$ is a process defined on an extension $(\Omega ,%
\mathcal{\tilde{F}},(\mathcal{\tilde{F})}_{t\geq 0},\mathbb{P})$ of $(\Omega
,\mathcal{F},\left( \mathcal{F}\right) _{t\geq 0},\mathbb{P})$, which
conditionally on $\mathcal{F}$ is a continuous Gaussian process, with
independent increments, and whose mean and variance are given respectively by%
\begin{eqnarray*}
\mathbb{\tilde{E}}[\tilde{X}_{2,t}|\mathcal{F]} &=&\frac{1}{2\sigma ^{2}}%
\lambda (\varphi _{p,\sigma })\sum_{j\geq 1}\sum_{k>j}\int_{0}^{t}\mathbb{I}%
_{\{\wedge _{l=j,k}Z_{l,s}>\vee _{l\neq j,k}Z_{l,s}\}}dL_{Z_{k\backslash
j},s}^{0} \\
\mathbb{\tilde{V}}[\tilde{X}_{2,t}|\mathcal{F]} &=&(m_{2p}-m_{p}^{2})\sigma
^{2p}t.
\end{eqnarray*}
\end{Prop}

\begin{Rk}
We only consider integers $p$ for technical reasons in the proof of the
central limit theorem, although it is expected that the asymptotic
convergence still holds for any $p\geq 1$.
\end{Rk}

We now consider the case where $H_{t}$ is not necessarily a constant
function and study the power variations of%
\begin{equation*}
\log \eta _{t}=\bigvee_{i=1}^{\infty }\left( \log R_{i}+\int_{0}^{t}H_{s}%
\text{d}W_{i,s}-\frac{1}{2}\int_{0}^{t}H_{s}^{2}\text{d}s\right)
=\bigvee_{i=1}^{\infty }\left( \log R_{i}+\int_{0}^{t}H_{s}\text{d}%
W_{i,s}\right) -\frac{1}{2}\int_{0}^{t}H_{s}^{2}\text{d}s
\end{equation*}%
where $(R_{i})_{i\geq 1}$ is a $\mathcal{F}_{0}$-Poisson point process, with
mean measure $dr/r^{2}$, $W_{1},W_{2},...$ are independent $\left( \mathcal{F%
}_{t}\right) $-Brownian motions, $H$ is a H\"{o}lder function in $C^{+}[0,1]$
with exponent $\alpha >1/2$ and satisfying $\int_{0}^{1}H_{s}^{2}ds<\infty $%
. For $i,j,k\geq 1$, let us also use the following notation%
\begin{equation*}
Z_{i,t}=\log R_{i}+\int_{0}^{t}H_{s}\text{d}W_{i,s},\quad Z_{k\backslash
j,t}=Z_{k,t}-Z_{j,t}.
\end{equation*}

\begin{Prop}
\label{TV_logY}As $n\rightarrow \infty $, we have, for any integer $p\geq 1$,%
\begin{equation*}
B\left( p,\log \eta \right) _{t}^{n}\overset{u.c.p}{\Longrightarrow }%
m_{p}\int_{0}^{t}H_{s}^{p}ds\quad \text{and}\quad \sqrt{n}\left( B\left(
p,\log \eta \right) _{t}^{n}-m_{p}\int_{0}^{t}H_{s}^{p}ds\right) \overset{%
\mathcal{L}-s}{\implies }\tilde{X}_{3,t},
\end{equation*}%
where $\tilde{X}_{3}$ is a process defined on an extension $(\Omega ,%
\mathcal{\tilde{F}},(\mathcal{\tilde{F})}_{t\geq 0},\mathbb{P})$ of $(\Omega
,\mathcal{F},\left( \mathcal{F}\right) _{t\geq 0},\mathbb{P})$, which
conditionally on $\mathcal{F}$ is a continuous Gaussian process, with
independent increments, and whose mean and variance are given respectively by%
\begin{eqnarray*}
\mathbb{\tilde{E}}[\tilde{X}_{3,t}|\mathcal{F]} &=&\frac{1}{2}\sum_{j\geq
1}\sum_{k>j}\int_{0}^{t}\frac{\lambda (\varphi _{p,H_{s}})}{H_{s}^{2}}%
\mathbb{I}_{\{\wedge _{l=j,k}Z_{l,s}>\vee _{l\neq
j,k}Z_{l,s}\}}dL_{Z_{k\backslash j},s}^{0} \\
\mathbb{\tilde{V}}[\tilde{X}_{3,t}|\mathcal{F]} &=&\left(
m_{2p}-m_{p}^{2}\right) \int_{0}^{t}H_{s}^{2p}ds.
\end{eqnarray*}
\end{Prop}

\section{Appendix}

$C$ is a constant that does not depend of $n$ but can vary from line to line.

\subsection{Proof of Proposition \protect\ref{TCLMaxB}}

We only prove the stable convergence in law of $\sqrt{n}\left( B\left(
p,W_{1}\vee W_{2}\right) _{t}^{n}-m_{p}t\right) $. We have%
\begin{equation*}
\sqrt{n}\left( B\left( p,W_{1}\vee W_{2}\right) _{t}^{n}-m_{p}t\right)
=\sum_{i=1}^{\left\lfloor nt\right\rfloor -1}\zeta
_{1,i}^{n}+\sum_{i=1}^{\left\lfloor nt\right\rfloor -1}\zeta
_{2,i}^{n}+m_{p}n^{1/2}\left( \frac{\left\lfloor nt\right\rfloor -1}{n}%
-t\right) ,
\end{equation*}%
where%
\begin{eqnarray*}
\zeta _{1,i}^{n} &=&n^{-1/2}(|\sqrt{n}\Delta _{i}^{n}W_{1}\mathbb{I}%
_{\{W_{2\backslash 1,(i-1)/n}<0\}}+\sqrt{n}\Delta _{i}^{n}W_{2}\mathbb{I}%
_{\{W_{2\backslash 1,(i-1)/n}>0\}}|^{p}-m_{p}) \\
\zeta _{2,i}^{n} &=&n^{-1/2}\Psi _{|\cdot |^{p}}(\sqrt{n}\Delta
_{i}^{n}W_{1},\sqrt{n}\Delta _{i}^{n}W_{2},\sqrt{n}W_{2\backslash
1,(i-1)/n}).
\end{eqnarray*}

Step 1) First it is clear that%
\begin{equation*}
n^{1/2}\left( \frac{\left\lfloor nt\right\rfloor -1}{n}-t\right) \overset{%
u.c.p}{\Longrightarrow }0.
\end{equation*}

Step 2) Second, by using usual arguments (see e.g. Chapter 5.2 in Jacod and
Protter (2011)), we have, as $n\rightarrow \infty $,%
\begin{equation*}
\sum_{i=1}^{\left\lfloor nt\right\rfloor -1}\zeta _{1,i}^{n}\overset{%
\mathcal{L}-s}{\implies }\tilde{X}_{0,t},
\end{equation*}%
where $\tilde{X}_{0}$ is a process defined on an extension $(\Omega ,%
\mathcal{\tilde{F}},(\mathcal{\tilde{F})}_{t\geq 0},\mathbb{P})$ of $(\Omega
,\mathcal{F},\left( \mathcal{F}\right) _{t\geq 0},\mathbb{P})$, which
conditionally on $\mathcal{F}$ is a continuous centered Gaussian process,
with independent increments, and whose variance is given by%
\begin{equation*}
\mathbb{\tilde{V}}[\tilde{X}_{0,t}|\mathcal{F]}=\left(
m_{2p}-m_{p}^{2}\right) t.
\end{equation*}

Step 3) Third, let us prove that 
\begin{equation*}
\sum_{i=1}^{\left\lfloor nt\right\rfloor -1}\zeta _{1,i}^{n}\overset{u.c.p}{%
\Longrightarrow }\frac{1}{2}\lambda (\varphi _{p})L_{W_{2\backslash
1},t}^{0}.
\end{equation*}%
We have%
\begin{eqnarray*}
\mathbb{E}\left[ \left. \zeta _{1,i}^{n}\right\vert \mathcal{F}_{(i-1)/n}%
\right] &=&n^{-1/2}\varphi _{p}\left( \sqrt{n}W_{2\backslash
1,(i-1)/n}\right) \\
\mathbb{E}\left[ \left. (\zeta _{1,i}^{n})^{2}\right\vert \mathcal{F}%
_{(i-1)/n}\right] &=&n^{-1}\varphi _{p}^{(2)}\left( \sqrt{n}W_{2\backslash
1,(i-1)/n}\right)
\end{eqnarray*}%
where%
\begin{equation*}
\varphi _{p}^{(2)}\left( w\right) =\int_{\mathbb{R}^{2}}\Psi _{|\cdot
|^{p}}^{2}(x,y,w)\frac{1}{2\pi }e^{-(x^{2}+y^{2})/2}\text{d}x\text{d}y.
\end{equation*}%
Note that $\int |\varphi _{2}\left( w\right) |$d$w<\infty $. It follows by
Theorem 1.1 in Jacod (1998) that%
\begin{eqnarray*}
&&\sum_{i=1}^{\left\lfloor nt\right\rfloor -1}\mathbb{E}\left[ \left. \zeta
_{1,i}^{n}\right\vert \mathcal{F}_{(i-1)/n}\right] \overset{u.c.p}{%
\Longrightarrow }\frac{1}{2}\lambda (\varphi _{p})L_{W_{2\backslash 1},t}^{0}
\\
&&\sqrt{n}\sum_{i=1}^{\left\lfloor nt\right\rfloor -1}\mathbb{E}\left[
\left. (\zeta _{1,i}^{n})^{2}\right\vert \mathcal{F}_{(i-1)/n}\right] 
\overset{u.c.p}{\Longrightarrow }\frac{1}{2}\lambda (\varphi
_{p}^{(2)})L_{W_{2\backslash 1},t}^{0}.
\end{eqnarray*}%
Therefore%
\begin{equation*}
\sum_{i=1}^{\left\lfloor nt\right\rfloor -1}\mathbb{E}\left[ \left. (\zeta
_{1,i}^{n})^{2}\right\vert \mathcal{F}_{(i-1)/n}\right] \overset{u.c.p}{%
\Longrightarrow }0,
\end{equation*}%
and the result follows by using Lemma 2.2.12 in Jacod and Protter (2011).

Step 4) Use steps 1), 2) and 3), define $\tilde{X}_{1,t}=\lambda (\varphi
_{p})L_{W_{2\backslash 1},t}^{0}/2+\tilde{X}_{0,t}$ and take into account
the properties of the stable convergence in law to conclude.

\subsection{Proof of Proposition \protect\ref{PropIncMAS}}

1) We first characterize the conditional distribution of $U_{i}^{n}$ given $%
\eta _{(i-1)/n}=\eta $. By Proposition 4.1 in Dombry and Eyi Menko (2013),
the conditional distribution $\eta _{i/n}|\eta _{(i-1)/n}=\eta $ is given by%
\begin{eqnarray*}
&&\Pr \left( \eta _{i/n}\leq z|\eta _{(i-1)/n}=\eta \right) \\
&=&\exp \left( -\mathbb{E}\left[ \left( \frac{V_{i/n}}{z}-\frac{V_{(i-1)/n}}{%
\eta }\right) _{+}\right] \right) \mathbb{E}\left[ \mathbb{I}\left( \frac{%
V_{i/n}}{z}<\frac{V_{(i-1)/n}}{\eta }\right) V_{(i-1)/n}\right] ,
\end{eqnarray*}%
and, by stationarity of $\eta $, it can be rewritten is the following way%
\begin{equation*}
\Pr \left( \eta _{i/n}\leq z|\eta _{(i-1)/n}=\eta \right) =\exp \left( -%
\mathbb{E}\left[ \left( \frac{V_{1/n}}{z}-\frac{1}{\eta }\right) _{+}\right]
\right) \mathbb{E}\left[ \mathbb{I}\left( \frac{V_{1/n}}{z}<\frac{1}{\eta }%
\right) \right] .
\end{equation*}%
We have%
\begin{eqnarray*}
\mathbb{E}\left[ \mathbb{I}\left( V_{1/n}<\frac{z}{\eta }\right) \right]
&=&\Pr \left( V_{1/n}<\frac{z}{\eta }\right) \\
&=&\Pr \left( W_{1/n}<\sigma ^{-1}\ln \left( \frac{z}{\eta }\right) +\frac{%
\sigma }{2n}\right) \\
&=&\Phi \left( \frac{1}{\sigma n^{-1/2}}\ln \left( \frac{z}{\eta }\right) +%
\frac{1}{2}\sigma n^{-1/2}\right) ,
\end{eqnarray*}%
and%
\begin{equation*}
\mathbb{E}\left[ \left( \frac{V_{1/n}}{z}-\frac{1}{\eta }\right) _{+}\right]
=\frac{1}{z}\Phi \left( -\frac{1}{\sigma n^{-1/2}}\ln \left( \frac{z}{\eta }%
\right) +\frac{\sigma }{2\sqrt{n}}\right) -\frac{1}{\eta }\Phi \left( -\frac{%
1}{\sigma n^{-1/2}}\ln \left( \frac{z}{\eta }\right) -\frac{\sigma }{2\sqrt{n%
}}\right) .
\end{equation*}%
Since%
\begin{equation*}
\Pr (\left. U_{i}^{n}\leq u\right\vert \eta _{(i-1)/n}=\eta )=\Pr \left(
\eta _{i/n}\leq \eta e^{u\sigma n^{-1/2}}|\eta _{(i-1)/n}=\eta \right) ,
\end{equation*}%
we deduce that%
\begin{eqnarray*}
&&\Pr \left( \left. U_{i}^{n}\leq u\right\vert \eta _{(i-1)/n}=\eta \right)
\\
&=&\exp \left( -\frac{1}{\eta }\left[ e^{-\sigma u/\sqrt{n}}\Phi \left( -u+%
\frac{\sigma }{2\sqrt{n}}\right) -\Phi \left( -u-\frac{\sigma }{2\sqrt{n}}%
\right) \right] \right) \Phi \left( u+\frac{\sigma }{2\sqrt{n}}\right) .
\end{eqnarray*}

2) We have 
\begin{equation*}
\Pr \left( U_{i}^{n}\leq u\right) =\mathbb{E}_{\eta _{(i-1)/n}}\left[ \Pr
(\left. U_{i}^{n}\leq u\right\vert \eta _{(i-1)/n})\right] ,
\end{equation*}%
and since $\eta _{(i-1)/n}^{-1}$ has a standard Exponential distribution, we
derive that%
\begin{equation*}
\Pr \left( U_{i}^{n}\leq u\right) =\frac{\Phi \left( u+\sigma /(2\sqrt{n}%
)\right) }{\Phi \left( u+\sigma /(2\sqrt{n})\right) +e^{-\sigma u/\sqrt{n}%
}\left( 1-\Phi \left( u-\sigma /(2\sqrt{n})\right) \right) }.
\end{equation*}

3) We have for any $p\geq 1$%
\begin{eqnarray*}
\mathbb{E}\left[ |U_{i}^{n}|^{p}\right] &=&p\int_{-\infty }^{0}(-u)^{p-1}\Pr
\left( U_{i}^{n}\leq u\right) du+p\int_{0}^{\infty }u^{p-1}\Pr \left(
U_{i}^{n}>u\right) du \\
&=&2p\int_{0}^{\infty }u^{p-1}\frac{e^{-\sigma u/\sqrt{n}}\left( 1-\Phi
\left( u-\sigma /(2\sqrt{n})\right) \right) }{\Phi \left( u+\sigma /(2\sqrt{n%
})\right) +e^{-\sigma u/\sqrt{n}}\left( 1-\Phi \left( u-\sigma /(2\sqrt{n}%
)\right) \right) }du.
\end{eqnarray*}%
Now let us define%
\begin{equation*}
V(n,u)=\frac{e^{-\sigma n^{-1/2}u}\bar{\Phi}\left( u-\sigma
n^{-1/2}/2\right) }{\bar{\Phi}\left( u\right) }\quad \text{and}\quad U(n,u)=%
\frac{\Phi \left( u+\sigma n^{-1/2}/2\right) }{\Phi \left( u\right) },
\end{equation*}%
and note that%
\begin{equation*}
\mathbb{E}\left[ |U_{i}^{n}|^{p}\right] =2p\int_{0}^{\infty }u^{p-1}\frac{%
\bar{\Phi}\left( u\right) V(n,u)}{\Phi \left( u\right) U(n,u)+\bar{\Phi}%
\left( u\right) V(nt,u)}du.
\end{equation*}

i) There exists $u_{c}\in \lbrack u,u+\sigma n^{-1/2}/2]$ such that 
\begin{equation*}
\Phi \left( u+\sigma /(2\sqrt{n})\right) =\Phi \left( u\right) +\sigma \frac{%
1}{\sqrt{n}}\varphi (u)/2+\frac{1}{8}\sigma ^{2}\frac{1}{n}\varphi ^{\prime
}(u_{c}).
\end{equation*}%
Therefore%
\begin{eqnarray*}
\left\vert \frac{\Phi \left( u+\sigma /(2\sqrt{n})\right) }{\Phi \left(
u\right) }-1-\sigma \frac{1}{\sqrt{n}}\frac{\varphi (u)}{2\Phi \left(
u\right) }\right\vert &\leq &\frac{1}{8}\sigma ^{2}\frac{1}{n}\frac{1}{\Phi
\left( u\right) }\sup_{u_{c}\in \lbrack u,u+\sigma n^{-1/2}/2]}u_{c}\varphi
(u_{c}) \\
&\leq &\frac{1}{8}\sigma ^{2}\frac{1}{n}\frac{1}{\Phi \left( u\right) }%
\left( u+\sigma n^{-1/2}/2\right) \varphi (u) \\
&\leq &\frac{1}{8}\sigma ^{2}\frac{1}{n}\sup_{u\geq 0}\frac{u\varphi (u)}{%
\Phi \left( u\right) }+\frac{1}{16}\sigma ^{3}\frac{1}{n^{3/2}}\sup_{u\geq 0}%
\frac{\varphi (u)}{\Phi \left( u\right) } \\
&\leq &O_{u\in \lbrack 0,\infty )}(\frac{1}{n}),
\end{eqnarray*}%
and we can deduce that 
\begin{equation*}
U(n,u)=1+\sigma \frac{1}{\sqrt{n}}\frac{\varphi (u)}{2\Phi \left( u\right) }+%
\tilde{U}(n,u)
\end{equation*}%
where%
\begin{equation*}
\sup_{u\in \lbrack 0,\infty )}\left\vert \tilde{U}(n,u)\right\vert \leq 
\frac{C}{n}.
\end{equation*}

ii) If $u\in \lbrack 0,n^{\beta })$ with $\beta <1/2$, then%
\begin{equation*}
e^{-\sigma n^{-1/2}u}=1-\sigma n^{-1/2}u+o_{u\in \lbrack 0,n^{\beta })}(%
\frac{1}{n}).
\end{equation*}%
There exists $u_{c}\in \lbrack u-\sigma n^{-1/2}/2,u]$ such that 
\begin{equation*}
\bar{\Phi}\left( u-\sigma n^{-1/2}/2\right) =\bar{\Phi}\left( u\right)
+\sigma n^{-1/2}\varphi (u)/2-\frac{1}{8}\sigma ^{2}\frac{1}{n}\varphi
^{\prime }(u_{c}).
\end{equation*}%
Therefore, if moreover $\beta <1/4$, we have 
\begin{eqnarray*}
\left\vert \frac{\bar{\Phi}\left( u-\sigma n^{-1/2}/2\right) }{\bar{\Phi}%
\left( u\right) }-1-\sigma n^{-1/2}\frac{\varphi (u)}{2\bar{\Phi}\left(
u\right) }\right\vert &\leq &\frac{1}{8}\sigma ^{2}\frac{1}{n}\frac{1}{\bar{%
\Phi}\left( u\right) }\sup_{u_{c}\in \lbrack u-\sigma
n^{-1/2}/2,u]}u_{c}\varphi (u_{c}) \\
&\leq &\frac{1}{8}\sigma ^{2}\frac{1}{n}\frac{1}{\bar{\Phi}\left( u\right) }%
u\varphi (u-\sigma n^{-1/2}/2) \\
&\leq &\frac{1}{8}\sigma ^{2}\frac{1}{n}\sup_{u\in \lbrack 0,n^{\beta })}%
\frac{u\varphi (u)}{\bar{\Phi}\left( u\right) }e^{\sigma n^{-1/2}u} \\
&\leq &\frac{1}{8}\sigma ^{2}\frac{1}{n}\left( n^{2\beta }e^{\sigma
n^{-(1/2-\beta )}}\right) \\
&\leq &o_{u\in \lbrack 0,n^{\beta })}(n^{-1/2}).
\end{eqnarray*}%
If $u\in \lbrack 0,n^{\beta })$ with $\beta <1/4$, we deduce that%
\begin{equation*}
V(n,u)=1+\sigma \frac{1}{\sqrt{n}}\left( \frac{\varphi (u)}{2\bar{\Phi}%
\left( u\right) }-u\right) +\tilde{V}(n,u),
\end{equation*}%
where%
\begin{equation*}
\sup_{u\in \lbrack 0,n^{\beta })}\left\vert \tilde{V}(n,u)\right\vert
=o(n^{-1/2}).
\end{equation*}%
Moreover, as $u\rightarrow \infty $ and $n$ is fixed, we have%
\begin{equation*}
V(n,u)=\frac{e^{-\sigma n^{-1/2}u}\bar{\Phi}\left( u-\sigma
n^{-1/2}/2\right) }{\bar{\Phi}\left( u\right) }=\exp \left( -\sigma \frac{u}{%
2\sqrt{n}}-\frac{1}{8}\frac{\sigma ^{2}}{n}+o(1)\right) .
\end{equation*}

iii) Now we have 
\begin{eqnarray*}
&&\int_{0}^{\infty }u^{p-1}\frac{\bar{\Phi}\left( u\right) V(n,u)}{\Phi
\left( u\right) U(n,u)+\bar{\Phi}\left( u\right) V(nt,u)}du \\
&=&\int_{0}^{n^{\beta }}u^{p-1}\frac{\bar{\Phi}\left( u\right) V(n,u)}{\Phi
\left( u\right) U(n,u)+\bar{\Phi}\left( u\right) V(n,u)}du+\int_{n^{\beta
}}^{\infty }u^{p-1}\frac{\bar{\Phi}\left( u\right) V(n,u)}{\Phi \left(
u\right) U(n,u)+\bar{\Phi}\left( u\right) V(n,u)}du.
\end{eqnarray*}%
Then, for $\beta <1/4$, 
\begin{eqnarray*}
&&\int_{0}^{n^{\beta }}u^{p-1}\frac{\bar{\Phi}\left( u\right) V(n,u)}{\Phi
\left( u\right) U(n,u)+\bar{\Phi}\left( u\right) V(n,u)}du \\
&=&\int_{0}^{n^{\beta }}u^{p-1}\frac{\bar{\Phi}\left( u\right) \left(
1+\sigma n^{-1/2}\left( \frac{\varphi (u)}{2\bar{\Phi}\left( u\right) }%
-u\right) +\tilde{V}(n,u)\right) }{\Phi \left( u\right) \left( 1+\sigma
n^{-1/2}\frac{\varphi (u)}{2\Phi \left( u\right) }+\tilde{U}(n,u)\right) +%
\bar{\Phi}\left( u\right) \left( 1+\sigma n^{-1/2}\left( \frac{\varphi (u)}{2%
\bar{\Phi}\left( u\right) }-u\right) +\tilde{V}(n,u)\right) }du \\
&=&\int_{0}^{n^{\beta }}u^{p-1}\frac{\bar{\Phi}\left( u\right) \left(
1+\sigma n^{-1/2}\left( \frac{\varphi (u)}{2\bar{\Phi}\left( u\right) }%
-u\right) +\tilde{V}(n,u)\right) }{1+\sigma n^{-1/2}(\varphi (u)-u\bar{\Phi}%
\left( u\right) )+\Phi \left( u\right) \tilde{U}(n,u)+\bar{\Phi}\left(
u\right) \tilde{V}(n,u)}du \\
&=&\int_{0}^{n^{\beta }}u^{p-1}\bar{\Phi}\left( u\right) du+\sigma
n^{-1/2}\int_{0}^{n^{\beta }}u^{p-1}\varphi (u)(1/2-\bar{\Phi}\left(
u\right) -u\bar{\Phi}\left( u\right) \Phi \left( u\right) /\varphi
(u))du+W_{n}
\end{eqnarray*}%
where $W_{n}=o(n^{-1/2})$. Moreover%
\begin{equation*}
\int_{n^{\beta }}^{\infty }u^{p-1}\frac{\bar{\Phi}\left( u\right) V(n,u)}{%
\Phi \left( u\right) U(n,u)+\bar{\Phi}\left( u\right) V(n,u)}du\leq
\int_{n^{\beta }}^{\infty }u^{p-1}\frac{\bar{\Phi}\left( u\right) }{\Phi
\left( u\right) }du\leq C\int_{n^{\beta }}^{\infty }u^{p-1}\bar{\Phi}\left(
u\right) du).
\end{equation*}%
Finally, we deduce that 
\begin{equation*}
\lim_{n\rightarrow \infty }\sqrt{n}\left( \mathbb{E}\left[ |U_{i}^{n}|^{p}%
\right] -2p\int_{0}^{\infty }u^{p-1}\bar{\Phi}\left( u\right) du\right)
=2p\int_{0}^{\infty }u^{p-1}\varphi (u)[1/2-\bar{\Phi}\left( u\right) -u\bar{%
\Phi}\left( u\right) \Phi \left( u\right) /\varphi (u)]du.
\end{equation*}

\subsection{Proof of Proposition \protect\ref{TCLSimpMS}}

We only prove the stable convergence in law of $\sqrt{n}\left( B\left(
p,\log \eta \right) _{t}^{n}-m_{p}\sigma ^{p}t\right) $. Recall that 
\begin{equation*}
\log \eta _{t}=\bigvee_{i=1}^{\infty }\left( \log R_{i}+\sigma W_{i,t}-\frac{%
1}{2}\sigma ^{2}t\right) =\bigvee_{i=1}^{\infty }Z_{i,t}-\frac{1}{2}\sigma
^{2}t,
\end{equation*}%
and that%
\begin{eqnarray*}
&&\left\vert \sqrt{n}\Delta _{i}^{n}\log \eta +\frac{1}{2\sqrt{n}}\sigma
^{2}\right\vert ^{p} \\
&=&\left\vert \sigma \sum_{j\geq 1}\sqrt{n}\Delta _{i}^{n}W_{j}\mathbb{I}%
_{\{\vee _{l\geq 1}Z_{l,(i-1)/n}=Z_{j,(i-1)/n}\}}\right\vert ^{p} \\
&&+\sum_{j\geq 1}\sum_{k>j}\Psi _{\left\vert \cdot \right\vert ^{p},\sigma }(%
\sqrt{n}\Delta _{i}^{n}W_{j},\sqrt{n}\Delta _{i}^{n}W_{k},\sqrt{n}%
Z_{k\backslash j,(i-1)/n}\mathbb{I}_{\{\vee _{l\neq j,k}Z_{l,(i-1)/n}\leq
\wedge _{l=j,k}Z_{l,(i-1)/n}\}}+H_{\left\vert \cdot \right\vert ^{p},i}^{n},
\end{eqnarray*}%
where 
\begin{eqnarray*}
H_{\left\vert \cdot \right\vert ^{p},i}^{n} &=&\sum_{j\geq 1}\sum_{k\neq
j}\Psi _{\left\vert \cdot \right\vert ^{p},\sigma }^{<}(\sqrt{n}\Delta
_{i}^{n}W_{j},\sqrt{n}\Delta _{i}^{n}W_{k},\sqrt{n}Z_{k\backslash
j,(i-1)/n})\times \\
&&\dprod\limits_{l\neq j,k}\mathbb{I}_{\{Z_{l\backslash j,(i-1)/n}\leq
0\}}\times \left( \dprod\limits_{l\neq k,j}\mathbb{I}_{\{Z_{l\backslash
k,i/n}\leq 0\}}-\dprod\limits_{l\neq k,j}\mathbb{I}_{\{Z_{l\backslash
k,(i-1)/n}\leq 0\}}\right) .
\end{eqnarray*}%
Note that, for $a,u\in \mathbb{R}$,%
\begin{equation*}
|a+u|=|a|+sign(a)u-2(|a|+sign(a)u)\mathbb{I}_{\{(a+u)a<0\}},
\end{equation*}%
and, for any integer $p\geq 1$,%
\begin{eqnarray*}
|a+u|^{p} &=&|a|^{p}+p\left[ sign(a)u-2(|a|+sign(a)u)\mathbb{I}%
_{\{(a+u)a<0\}}\right] |a|^{p-1} \\
&&+\sum_{k=2}^{p}C_{p}^{k}\left[ sign(a)u-2(|a|+sign(a)u)\mathbb{I}%
_{\{(a+u)a<0\}}\right] ^{k}|a|^{p-k}
\end{eqnarray*}%
where $C_{p}^{k}=p!/(k!(p-k)!)$ is the binomial coefficient of order $k$ and 
$p$.

Therefore we have%
\begin{eqnarray*}
B\left( p,\log \eta \right) _{t}^{n} &=&\frac{1}{n}\sum_{i=1}^{\left\lfloor
nt\right\rfloor -1}\left\vert \sqrt{n}\Delta _{i}^{n}\log \eta \right\vert
^{p} \\
&=&\frac{1}{n}\sum_{i=1}^{\left\lfloor nt\right\rfloor -1}\left\vert \sqrt{n}%
\Delta _{i}^{n}\log \eta +\frac{1}{2\sqrt{n}}\sigma ^{2}\right\vert
^{p}+A_{t}^{n}+B_{t}^{n}+C_{t}^{n},
\end{eqnarray*}%
where 
\begin{eqnarray*}
A_{t}^{n} &=&-\frac{1}{2n^{3/2}}\sigma ^{2}p\sum_{i=1}^{\left\lfloor
nt\right\rfloor -1}sign\left( \Delta _{i}^{n}\log \eta \right) \left\vert 
\sqrt{n}\Delta _{i}^{n}\log \eta \right\vert ^{p-1} \\
B_{t}^{n} &=&\frac{2}{n}p\sum_{i=1}^{\left\lfloor nt\right\rfloor -1}\left( |%
\sqrt{n}\Delta _{i}^{n}\log \eta |+sign(\Delta _{i}^{n}\log \eta )\frac{1}{2%
\sqrt{n}}\sigma ^{2}\right) \left\vert \sqrt{n}\Delta _{i}^{n}\log \eta
\right\vert ^{p-1}\mathbb{I}_{\{(\sqrt{n}\Delta _{i}^{n}\log \eta
+n^{-1/2}\sigma ^{2}/2)\sqrt{n}\Delta _{i}^{n}\log \eta <0\}} \\
C_{t}^{n} &=&\frac{1}{n}\sum_{i=1}^{\left\lfloor nt\right\rfloor
-1}\sum_{i=1}^{\left\lfloor nt\right\rfloor -1}C_{p}^{k}|\sqrt{n}\Delta
_{i}^{n}\log \eta |^{p-k}\times \\
&&\left[ sign\left( \Delta _{i}^{n}\log \eta \right) \frac{1}{2\sqrt{n}}%
\sigma ^{2}-2(|\sqrt{n}\Delta _{i}^{n}\log \eta |+sign(\Delta _{i}^{n}\log
\eta )\frac{1}{2\sqrt{n}}\sigma ^{2})\mathbb{I}_{\{(\sqrt{n}\Delta
_{i}^{n}\log \eta +n^{-1/2}\sigma ^{2}/2)\sqrt{n}\Delta _{i}^{n}\log \eta
<0\}}\right] ^{k}
\end{eqnarray*}%
It follows that%
\begin{eqnarray*}
&&\sqrt{n}\left( B\left( p,\log \eta \right) _{t}^{n}-m_{p}\sigma
^{p}t\right) \\
&=&\sum_{i=1}^{\left\lfloor nt\right\rfloor -1}\zeta
_{1,i}^{n}+\sum_{i=1}^{\left\lfloor nt\right\rfloor -1}\zeta _{2,i}^{n}+%
\frac{1}{\sqrt{n}}\sum_{i=1}^{\left\lfloor nt\right\rfloor -1}H_{\left\vert
\cdot \right\vert ^{p},i}^{n} \\
&&+\sqrt{n}A_{t}^{n}+\sqrt{n}B_{t}^{n}+\sqrt{n}C_{t}^{n}+m_{p}\sigma
^{p}n^{1/2}\left( \frac{\left\lfloor nt\right\rfloor -1}{n}-t\right)
\end{eqnarray*}%
where%
\begin{eqnarray*}
\zeta _{1,i}^{n} &=&n^{-1/2}\sigma (|\sum_{j\geq 1}\sqrt{n}\Delta
_{i}^{n}W_{j}\mathbb{I}_{\{\vee _{l\geq
1}Z_{l,(i-1)/n}=Z_{j,(i-1)/n}\}}|^{p}-m_{p}) \\
\zeta _{2,i}^{n} &=&n^{-1/2}\sum_{j\geq 1}\sum_{k>j}\Psi _{\left\vert \cdot
\right\vert ^{p},\sigma }(\sqrt{n}\Delta _{i}^{n}W_{j},\sqrt{n}\Delta
_{i}^{n}W_{k},\sqrt{n}Z_{k\backslash j,(i-1)/n}\mathbb{I}_{\{\vee _{l\neq
j,k}Z_{l,(i-1)/n}\leq \wedge _{l=j,k}Z_{l,(i-1)/n}\}}.
\end{eqnarray*}

Step 1) First it is clear that%
\begin{equation*}
n^{1/2}\left( \frac{\left\lfloor nt\right\rfloor -1}{n}-t\right) \overset{%
u.c.p}{\Longrightarrow }0.
\end{equation*}

Step 2) We have%
\begin{equation*}
\sqrt{n}A_{t}^{n}=\sum_{i=1}^{\left\lfloor nt\right\rfloor -1}\xi _{i}^{n}
\end{equation*}%
where%
\begin{equation*}
\xi _{i}^{n}=-\frac{\sigma ^{2}}{2n}p\times sign\left( \Delta _{i}^{n}\log
\eta \right) \left\vert \sqrt{n}\Delta _{i}^{n}\log \eta \right\vert ^{p-1}.
\end{equation*}%
Note that, by Proposition \ref{PropIncMAS}, we have for $p=1$%
\begin{eqnarray*}
\mathbb{E}\left[ \left. \xi _{i}^{n}\right\vert \eta _{(i-1)/n}\right] &=&-%
\frac{\sigma ^{2}}{2n}\left[ \Pr (\left. U_{i}^{n}>0\right\vert \eta
_{(i-1)/n})-\Pr (\left. U_{i}^{n}\leq 0\right\vert \eta _{(i-1)/n})\right] \\
&=&-\frac{\sigma ^{2}}{4n}\left[ \frac{1}{2}-\Pr (\left. U_{i}^{n}\leq
0\right\vert \eta _{(i-1)/n})\right] \\
&=&-\frac{\sigma ^{2}}{4n}\left[ \frac{1}{2}-\exp \left( -\frac{1}{\eta
_{(i-1)/n}}\left[ \Phi \left( \frac{\sigma }{2\sqrt{n}}\right) -\Phi \left( -%
\frac{\sigma }{2\sqrt{n}}\right) \right] \right) \Phi \left( \frac{\sigma }{2%
\sqrt{n}}\right) \right] \\
&=&-\frac{\sigma ^{2}}{8n}\left[ \frac{1}{\eta _{(i-1)/n}}-1\right] \frac{%
\sigma }{\sqrt{n}}\frac{1}{\sqrt{2\pi }}\left( 1+o(1)\right) ,
\end{eqnarray*}%
and, for any integer larger than $1$, the same type of calculations leads to%
\begin{equation*}
\mathbb{E}\left[ \left. \xi _{i}^{n}\right\vert \eta _{(i-1)/n}\right] =%
\frac{C}{n^{3/2}}\left[ \frac{1}{\eta _{(i-1)/n}}-1\right] \left(
1+o(1)\right) .
\end{equation*}%
Therefore we have%
\begin{equation*}
\sum_{i=1}^{\left\lfloor nt\right\rfloor -1}\mathbb{E}\left[ \left. \xi
_{i}^{n}\right\vert \eta _{(i-1)/n}\right] \overset{u.c.p}{\Longrightarrow }%
0.
\end{equation*}%
Moreover we have%
\begin{equation*}
\mathbb{E}\left[ \left. (\xi _{i}^{n})^{2}\right\vert \eta _{(i-1)/n}\right]
=\frac{\sigma ^{4}}{4n^{2}}m_{2p-2}\left( 1+o(1)\right)
\end{equation*}%
and therefore%
\begin{equation*}
\sum_{i=1}^{\left\lfloor nt\right\rfloor -1}\mathbb{E}\left[ \left. (\xi
_{i}^{n})^{2}\right\vert \mathcal{F}_{(i-1)/n}^{\eta }\right] \overset{u.c.p}%
{\Longrightarrow }0.
\end{equation*}%
By Lemma 2.2.10 in Jacod and Protter (2011), we deduce that 
\begin{equation*}
\sqrt{n}A_{t}^{n}\overset{u.c.p}{\Longrightarrow }0.
\end{equation*}

Step 3) We have%
\begin{equation*}
\sqrt{n}B_{t}^{n}=\sum_{i=1}^{\left\lfloor nt\right\rfloor -1}\zeta _{i}^{n}
\end{equation*}%
with%
\begin{equation*}
\zeta _{i}^{n}=\frac{2}{\sqrt{n}}\left( |\sqrt{n}\Delta _{i}^{n}\log \eta
|+sign(\Delta _{i}^{n}\log \eta )\frac{1}{2\sqrt{n}}\sigma ^{2}\right)
\left\vert \sqrt{n}\Delta _{i}^{n}\log \eta \right\vert ^{p-1}\mathbb{I}_{\{(%
\sqrt{n}\Delta _{i}^{n}\log \eta +n^{-1/2}\sigma ^{2}/2)\sqrt{n}\Delta
_{i}^{n}\log \eta <0\}}.
\end{equation*}%
Note that, for $a,u\in \mathbb{R}$,%
\begin{equation*}
||a|+sign(a)u|\mathbb{I}_{\{(a+u)a<0\}}\leq 2|a|\mathbb{I}_{\{|a|<|u|\}}\leq
2|u|\mathbb{I}_{\{|a|<|u|\}}.
\end{equation*}%
Therefore we have%
\begin{eqnarray*}
\mathbb{E}\left[ \left. |\zeta _{i}^{n}|\right\vert \eta _{(i-1)/n}\right]
&\leq &\frac{4}{\sqrt{n}}\left( \frac{\sigma ^{2}}{2\sqrt{n}}\right) ^{p}%
\left[ \Pr \left( \left. |\sqrt{n}\Delta _{i}^{n}\log \eta |<\frac{\sigma
^{2}}{2\sqrt{n}}\right\vert \eta _{(i-1)/n}\right) \right] \\
&=&\frac{2^{2-p}\sigma ^{2p}}{n^{(p+1)/2}}\Pr \left( \left. |U_{i}^{n}|<%
\frac{\sigma }{2\sqrt{n}}\right\vert \mathcal{F}_{(i-1)/n}^{\eta }\right)
\end{eqnarray*}%
where $U_{i}^{n}=\sigma ^{-1}\sqrt{n}\Delta _{i}^{n}\log \eta $. By
Proposition \ref{PropIncMAS}, we derive that%
\begin{eqnarray*}
\mathbb{E}\left[ \left. |\zeta _{i}^{n}|\right\vert \eta _{(i-1)/n}\right]
&\leq &\frac{2^{2-p}\sigma ^{2p}}{n^{(p+1)/2}}\exp \left( -\frac{1}{\eta
_{(i-1)/n}}\left[ \frac{e^{-\sigma ^{2}/2n}}{2}-\Phi \left( -\frac{\sigma
^{2}}{n}\right) \right] \right) \Phi \left( \frac{\sigma ^{2}}{n}\right) \\
&&-\frac{2^{2-p}\sigma ^{2p}}{n^{(p+1)/2}}\exp \left( -\frac{1}{\eta
_{(i-1)/n}}\left[ e^{\sigma ^{2}/n}\Phi \left( \frac{\sigma ^{2}}{n}\right) -%
\frac{1}{2}\right] \right) \frac{1}{2} \\
&\leq &C\frac{1}{n^{(p+3)/2}},
\end{eqnarray*}%
for large $n$. By Lemma 2.2.10 in Jacod and Protter (2011), it follows that 
\begin{equation*}
\sqrt{n}B_{t}^{n}=\sum_{i=1}^{\left\lfloor nt\right\rfloor -1}\zeta _{i}^{n}%
\overset{u.c.p}{\Longrightarrow }0.
\end{equation*}

Step 4) By using the same type of arguments as in Steps 2) and 3), it is
easily seen that%
\begin{equation*}
\sqrt{n}C_{t}^{n}\overset{u.c.p}{\Longrightarrow }0.
\end{equation*}

Step 5) Let us prove that 
\begin{equation*}
\sum_{i=1}^{\left\lfloor nt\right\rfloor -1}\zeta _{2,i}^{n}\overset{u.c.p}{%
\Longrightarrow }\frac{1}{2\sigma ^{2}}\lambda (\varphi _{p,\sigma
})\sum_{j\geq 1}\sum_{k>j}\int_{0}^{t}\mathbb{I}_{\{\wedge
_{l=j,k}Z_{l,s}>\vee _{l\neq j,k}Z_{l,s}\}}dL_{Z_{k\backslash j},s}^{0}.
\end{equation*}%
Since%
\begin{equation*}
\sum_{j\geq 1}\sum_{k>j}\mathbb{I}_{\{\vee _{l\neq j,k}Z_{l,(i-1)/n}\leq
\wedge _{l=j,k}Z_{l,(i-1)/n}\}}=1,
\end{equation*}%
it is enough to prove that for some $(j,k)$ such that $j\geq 1$ and $k>j$ 
\begin{eqnarray*}
&&\frac{1}{\sqrt{n}}\sum_{i=1}^{\left\lfloor nt\right\rfloor -1}\Psi
_{\left\vert \cdot \right\vert ^{p},\sigma }(\sqrt{n}\Delta _{i}^{n}W_{j},%
\sqrt{n}\Delta _{i}^{n}W_{k},\sqrt{n}Z_{k\backslash j,(i-1)/n})\mathbb{I}%
_{\{\vee _{l\neq j,k}Z_{l,(i-1)/n}\leq \wedge _{l=j,k}Z_{l,(i-1)/n}\}} \\
&&\overset{u.c.p}{\Longrightarrow }\frac{1}{2\sigma ^{2}}\lambda (\varphi
_{p,\sigma })\int_{0}^{t}\mathbb{I}_{\{\wedge _{l=j,k}Z_{l,s}>\vee _{l\neq
j,k}Z_{l,s}\}}dL_{Z_{k\backslash j},s}^{0}.
\end{eqnarray*}%
First note that%
\begin{equation*}
\mathbb{E}\left[ \left. \Psi _{\left\vert \cdot \right\vert ^{p},\sigma }(%
\sqrt{n}\Delta _{i}^{n}W_{j},\sqrt{n}\Delta _{i}^{n}W_{k},\sqrt{n}%
Z_{k\backslash j,(i-1)/n})\right\vert \mathcal{F}_{(i-1)/n}\right] =\varphi
_{p,\sigma }\left( \sqrt{n}Z_{k\backslash j,(i-1)/n}\right) .
\end{equation*}%
We can deduce from a simple modification of Theorem 1.1 in Jacod (1998), that%
\begin{eqnarray*}
&&\frac{1}{\sqrt{n}}\sum_{i=1}^{\left\lfloor nt\right\rfloor -1}\mathbb{E}%
\left[ \left. \Psi _{|\cdot |^{p},\sigma }(\sqrt{n}\Delta _{i}^{n}W_{j},%
\sqrt{n}\Delta _{i}^{n}W_{k},\sqrt{n}Z_{k\backslash j,(i-1)/n})\right\vert 
\mathcal{F}_{(i-1)/n}\right] \mathbb{I}_{\{\vee _{l\neq
j,k}Z_{l,(i-1)/n}\leq \wedge _{l=j,k}Z_{l,(i-1)/n}\}} \\
&&\overset{u.c.p}{\Longrightarrow }\frac{1}{2\sigma ^{2}}\lambda (\varphi
_{p,\sigma })\int_{0}^{t}\mathbb{I}_{\{\wedge _{l=j,k}Z_{l,s}>\vee _{l\neq
j,k}Z_{l,s}\}}dL_{Z_{k\backslash j},s}^{0}.
\end{eqnarray*}%
Moreover%
\begin{eqnarray*}
&&\frac{1}{\sqrt{n}}\sum_{i=1}^{\left\lfloor nt\right\rfloor -1}\mathbb{E}%
\left[ \left. \Psi _{|\cdot |^{p},\sigma }^{2}(\sqrt{n}\Delta _{i}^{n}W_{j},%
\sqrt{n}\Delta _{i}^{n}W_{k},\sqrt{n}Z_{k\backslash j,(i-1)/n})\right\vert 
\mathcal{F}_{(i-1)/n}\right] \mathbb{I}_{\{\vee _{l\neq
j,k}Z_{l,(i-1)/n}\leq \wedge _{l=j,k}Z_{l,(i-1)/n}\}} \\
&&\overset{u.c.p}{\Longrightarrow }\frac{1}{2\sigma ^{2}}\lambda (\varphi
_{p,\sigma }^{(2)})\int_{0}^{t}\mathbb{I}_{\{\wedge _{l=j,k}Z_{l,s}>\vee
_{l\neq j,k}Z_{l,s}\}}dL_{Z_{k\backslash j},s}^{0}
\end{eqnarray*}%
where%
\begin{equation*}
\varphi _{p,\sigma }^{(2)}\left( w\right) =\int_{\mathbb{R}^{2}}\Psi
_{|\cdot |^{p},\sigma }^{2}(x,y,w)\frac{1}{2\pi }e^{-(x^{2}+y^{2})/2}\text{d}%
x\text{d}y.
\end{equation*}%
The conclusion follows by using Lemma 2.2.12 in Jacod and Protter (2011).

Step 6) Let us prove that 
\begin{equation*}
\frac{1}{\sqrt{n}}\sum_{i=1}^{\left\lfloor nt\right\rfloor -1}H_{\left\vert
\cdot \right\vert ^{p},i}^{n}\overset{u.c.p}{\Longrightarrow }0.
\end{equation*}%
In the same way as in Step 3), it is enough to prove that for some $(j,k)$
such that $j\geq 1$, $k\geq 1$ and $k\neq j$,%
\begin{equation*}
\frac{1}{\sqrt{n}}\sum_{i=1}^{\left\lfloor nt\right\rfloor -1}\Psi
_{\left\vert \cdot \right\vert ^{p},\sigma }^{<}(\sqrt{n}\Delta
_{i}^{n}W_{j},\sqrt{n}\Delta _{i}^{n}W_{k},\sqrt{n}Z_{k\backslash
j,(i-1)/n})\times I_{i}^{(j,k),n}\overset{u.c.p}{\Longrightarrow }0,
\end{equation*}%
where%
\begin{equation*}
I_{i}^{(j,k),n}=\dprod\limits_{l\neq j,k}\mathbb{I}_{\{Z_{l\backslash
j,(i-1)/n}\leq 0\}}\times \left( \dprod\limits_{l\neq k,j}\mathbb{I}%
_{\{Z_{l\backslash k,i/n}\leq 0\}}-\dprod\limits_{l\neq k,j}\mathbb{I}%
_{\{Z_{l\backslash k,(i-1)/n}\leq 0\}}\right) .
\end{equation*}%
Let 
\begin{equation*}
Y_{k,j,t}=R_{k}+\sigma W_{k,t}-\vee _{l\neq k,j}\left( R_{l}+\sigma
W_{l,t}\right) .
\end{equation*}%
We have%
\begin{equation*}
\dprod\limits_{l\neq k,j}\mathbb{I}_{\{Z_{l\backslash k,i/n}\leq
0\}}-\dprod\limits_{l\neq k,j}\mathbb{I}_{\{Z_{l\backslash k,(i-1)/n}\leq
0\}}=\mathbb{I}_{\{Y_{k,j,(i-1)/n}>0,Y_{k,j,i/n}\leq 0\}}-\mathbb{I}%
_{\{Y_{k,j,(i-1)/n}\leq 0,Y_{k,j,i/n}>0\}}.
\end{equation*}%
Note that%
\begin{eqnarray*}
&&\Psi _{\left\vert \cdot \right\vert ^{p},\sigma }^{<}(\sqrt{n}\Delta
_{i}^{n}W_{j},\sqrt{n}\Delta _{i}^{n}W_{k},\sqrt{n}Z_{k\backslash j,(i-1)/n})
\\
&=&(|\sigma \sqrt{n}\Delta _{i}^{n}W_{k}+\sqrt{n}Z_{k\backslash
j,(i-1)/n}|^{p}-|\sigma \sqrt{n}\Delta _{i}^{n}W_{j}|^{p})\mathbb{I}%
_{\{\sigma \sqrt{n}\Delta _{i}^{n}W_{j}-\sigma \sqrt{n}\Delta
_{i}^{n}W_{k}\leq \sqrt{n}Z_{k\backslash j,(i-1)/n}\leq 0\}}.
\end{eqnarray*}%
Therefore%
\begin{equation*}
\Psi _{\left\vert \cdot \right\vert ^{p},\sigma }^{<}(\sqrt{n}\Delta
_{i}^{n}W_{j},\sqrt{n}\Delta _{i}^{n}W_{k},\sqrt{n}Z_{k\backslash
j,(i-1)/n})\times I_{i}^{(j,k),n}
\end{equation*}%
could be different from $0$, if at least $Z_{k\backslash j,(i-1)/n}$ is
close to zero, and $Y_{k,j,(i-1)/n}$ is also close to zero, or equivalently $%
Z_{k,(i-1)/n}$ is close to $Z_{j,(i-1)/n}$ and to $\vee _{l\neq
k,j}Z_{l,(i-1)/n}$. It is well known that bi-dimensional diffusion processes
never revisit a point in the plane, and so they do not in particular have a
local time. As a consequence, it is derived that%
\begin{equation*}
\frac{1}{\sqrt{n}}\sum_{i=1}^{\left\lfloor nt\right\rfloor -1}\Psi
_{\left\vert \cdot \right\vert ^{p},\sigma }^{<}(\sqrt{n}\Delta
_{i}^{n}W_{j},\sqrt{n}\Delta _{i}^{n}W_{k},\sqrt{n}Z_{k\backslash
j,(i-1)/n})\times I_{i}^{(j,k),n}\overset{u.c.p}{\Longrightarrow }0.
\end{equation*}

Step 7) By using usual arguments (see e.g. Chapter 5.2 in Jacod and Protter
(2011)), we have%
\begin{equation*}
\sum_{i=1}^{\left\lfloor nt\right\rfloor -1}\zeta _{1,i}^{n}\overset{%
\mathcal{L}-s}{\implies }\sigma \tilde{X}_{0,t},
\end{equation*}%
where $\tilde{X}_{0}$ is a process defined on an extension $(\Omega ,%
\mathcal{\tilde{F}},(\mathcal{\tilde{F})}_{t\geq 0},\mathbb{P})$ of $(\Omega
,\mathcal{F},\left( \mathcal{F}\right) _{t\geq 0},\mathbb{P})$, which
conditionally on $\mathcal{F}$ is a continuous centered Gaussian martingale
with variance 
\begin{equation*}
\mathbb{\tilde{V}}[\tilde{X}_{0,t}|\mathcal{F]}=\left(
m_{2p}-m_{p}^{2}\right) t.
\end{equation*}

Step 8) Use steps from 1) to 7) and define 
\begin{equation*}
\tilde{X}_{2,t}=\frac{1}{2\sigma ^{2}}\lambda (\varphi _{p,\sigma
})\sum_{j\geq 1}\sum_{k>j}\int_{0}^{t}\mathbb{I}_{\{\wedge
_{l=j,k}Z_{l,s}>\vee _{l\neq j,k}Z_{l,s}\}}dL_{Z_{k\backslash
j},s}^{0}+\sigma \tilde{X}_{0,t}
\end{equation*}%
to conclude.

\subsection{Proof of Proposition \protect\ref{TV_logY}}

We only prove the stable convergence in law of $\sqrt{n}(B\left( p,\log \eta
\right) _{t}^{n}-m_{p}\int_{0}^{t}H_{s}^{p}ds)$. Recall that 
\begin{equation*}
\log \eta _{t}=\bigvee_{j=1}^{\infty }\left( \log
R_{j}+\int_{0}^{t}H_{s}dW_{j,s}-\frac{1}{2}\int_{0}^{t}H_{s}^{2}ds\right)
=\bigvee_{j=1}^{\infty }Z_{j,t}-\frac{1}{2}\int_{0}^{t}H_{s}^{2}ds.
\end{equation*}%
Let%
\begin{equation*}
H_{s}^{(i)}=\left\{ 
\begin{tabular}{ll}
$H_{s}$ & if $s\leq (i-1)/n$ \\ 
$H_{(i-1)/n}$ & if $s>(i-1)/n$%
\end{tabular}%
\right. ,
\end{equation*}%
and define%
\begin{equation*}
\log \eta _{t}^{(i)}=\dbigvee\limits_{j=1}^{\infty }Z_{j,t}^{(i)}-\frac{1}{2}%
\int_{0}^{t}(H_{s}^{(i)})^{2}ds
\end{equation*}%
with%
\begin{equation*}
Z_{j,t}^{(i)}=\log R_{j}+\int_{0}^{t}H_{s}^{(i)}dW_{j,s}.
\end{equation*}

We have%
\begin{equation*}
\Delta _{i}^{n}\log \eta =\Delta _{i}^{n}\log \eta ^{(i)}+[\Delta
_{i}^{n}\log \eta -\Delta _{i}^{n}\log \eta ^{(i)}]
\end{equation*}%
Now, note that, for $a,u\in \mathbb{R}$,%
\begin{equation*}
|a+u|=|a|+sign(a)u-2(|a|+sign(a)u)\mathbb{I}_{\{(a+u)a<0\}}.
\end{equation*}%
With%
\begin{equation*}
a_{i}^{n}=\Delta _{i}^{n}\log \eta ^{(i)},\quad u_{i}^{n}=[\Delta
_{i}^{n}\log \eta -\Delta _{i}^{n}\log \eta ^{(i)}],
\end{equation*}%
and%
\begin{equation*}
b_{i}^{n}=sign(a_{i}^{n})u_{i}^{n}-2(|a_{i}^{n}|+sign(a_{i}^{n})u_{i}^{n})%
\mathbb{I}_{\{(a_{i}^{n}+u_{i}^{n})a_{i}^{n}<0\}},
\end{equation*}%
we have%
\begin{equation*}
|\Delta _{i}^{n}\log \eta |^{p}=|\Delta _{i}^{n}\log \eta _{i}|^{p}+w_{i}^{n}
\end{equation*}%
where%
\begin{equation*}
w_{i}^{n}=\sum_{k=1}^{p}C_{p}^{k}(b_{i}^{n})^{k}|a_{i}^{n}|^{p-k}.
\end{equation*}%
Therefore%
\begin{eqnarray*}
&&\sqrt{n}(B\left( p,\log \eta \right) _{t}^{n}-m_{p}\int_{0}^{t}H_{s}^{p}ds)
\\
&=&\sum_{i=1}^{\left\lfloor nt\right\rfloor -1}\zeta
_{i}^{n}+\sum_{i=1}^{\left\lfloor nt\right\rfloor
-1}n^{(p-1)/2}w_{i}^{n}+\sum_{i=1}^{\left\lfloor nt\right\rfloor
-1}x_{i}^{n}+\sqrt{n}m_{p}\int_{(\left\lfloor nt\right\rfloor
/n}^{t}H_{s}^{p}ds,
\end{eqnarray*}%
where%
\begin{eqnarray*}
\zeta _{i}^{n} &=&n^{-1/2}\left( |n^{1/2}\Delta _{i}^{n}\log \eta
^{(i)}|^{p}-m_{p}H_{(i-1)/n}^{p}\right) \\
x_{i}^{n} &=&-n^{1/2}m_{p}\int_{(i-1)/n}^{i/n}(H_{s}^{p}-H_{(i-1)/n}^{p})ds.
\end{eqnarray*}

Step 1) It is clear that%
\begin{equation*}
\sqrt{n}\int_{(\left\lfloor nt\right\rfloor /n}^{t}H_{s}^{p}ds\overset{u.c.p}%
{\Longrightarrow }0\text{.}
\end{equation*}

Step 2) Since $s\rightarrow H_{s}^{p}$ is H\"{o}lder with index $\alpha $,
we derive that%
\begin{equation*}
\int_{(i-1)/n}^{i/n}|H_{s}^{p}-H_{(i-1)/n}^{p}|ds\leq C\frac{1}{n^{1+\alpha }%
}
\end{equation*}%
and, as $n\rightarrow \infty $, we have%
\begin{equation*}
\sum_{i=1}^{\left\lfloor nt\right\rfloor -1}|x_{i}^{n}|\leq Cn^{1/2-\alpha
}\rightarrow 0\text{.}
\end{equation*}

Step 3)

i) We first establish that, for $q\geq 1$,%
\begin{equation*}
\mathbb{E}\left[ \left. |u_{i}^{n}|^{q}\right\vert \mathcal{F}_{(i-1)/n}%
\right] \leq C\frac{1}{n^{\alpha q+(2\wedge q)/2}}.
\end{equation*}%
Let $j_{t}$ be defined by%
\begin{equation*}
\log \eta _{t}=\dbigvee\limits_{j=1}^{\infty }Z_{j,t}=Z_{j_{t},t}.
\end{equation*}%
We have%
\begin{eqnarray*}
&&\left\vert \Delta _{i}^{n}\log \eta -\Delta _{i}^{n}\log \eta
_{i}\right\vert ^{q} \\
&\leq &\left\vert \int_{(i-1)/n}^{i/n}\left[ H_{s}-H_{(i-1)/n}\right]
dW_{j_{i/n},s}-\frac{1}{2}\int_{(i-1)/n}^{i/n}\left[
H_{s}^{2}-H_{(i-1)/n}^{2}\right] ds\right\vert ^{q} \\
&&+\left\vert \int_{(i-1)/n}^{i/n}\left[ H_{s}-H_{(i-1)/n}\right]
dW_{j_{(i-1)/n},s}-\frac{1}{2}\int_{(i-1)/n}^{i/n}\left[
H_{s}^{2}-H_{(i-1)/n}^{2}\right] ds\right\vert ^{q} \\
&\leq &C\left\vert \int_{(i-1)/n}^{i/n}\left[ H_{s}-H_{(i-1)/n}\right]
dW_{j_{i/n},s}\right\vert ^{q}+C\left\vert \int_{(i-1)/n}^{i/n}\left[
H_{s}-H_{(i-1)/n}\right] dW_{j_{(i-1)/n},s}\right\vert ^{q} \\
&&+C\int_{(i-1)/n}^{i/n}\left\vert H_{s}^{2}-H_{(i-1)/n}^{2}\right\vert
^{q}ds.
\end{eqnarray*}%
Then%
\begin{equation*}
\mathbb{E}\left[ \left. \left\vert \Delta _{i}^{n}\log \eta -\Delta
_{i}^{n}\log \eta _{i}\right\vert ^{q}\right\vert \mathcal{F}_{(i-1)/n}%
\right] \leq C\frac{1}{n^{\alpha q}}\frac{1}{n^{q/2}}+C\frac{1}{n^{\alpha q}}%
\frac{1}{n}\leq C\frac{1}{n^{\alpha q+(2\wedge q)/2}}.
\end{equation*}

ii) a) Study of $b_{i}^{n}$: 
\begin{equation*}
b_{i}^{n}=sign(a_{i}^{n})u_{i}^{n}-2(|a_{i}^{n}|+sign(a_{i}^{n})u_{i}^{n})%
\mathbb{I}_{\{(a_{i}^{n}+u_{i}^{n})a_{i}^{n}<0\}}:=b_{1,i}^{n}+b_{2,i}^{n}
\end{equation*}%
with%
\begin{eqnarray*}
b_{1,i}^{n} &=&sign(a_{i}^{n})u_{i}^{n} \\
b_{2,i}^{n} &=&-2(|a_{i}^{n}|+sign(a_{i}^{n})u_{i}^{n})\mathbb{I}%
_{\{(a_{i}^{n}+u_{i}^{n})a_{i}^{n}<0\}}.
\end{eqnarray*}%
For $q\geq 1$,%
\begin{equation*}
\mathbb{E}\left[ \left. |b_{i}^{n}|^{q}\right\vert \mathcal{F}_{(i-1)/n}%
\right] \leq C\left( \mathbb{E}\left[ \left. |b_{1,i}^{n}|^{q}\right\vert 
\mathcal{F}_{(i-1)/n}\right] +\mathbb{E}\left[ \left.
|b_{2,i}^{n}|^{q}\right\vert \mathcal{F}_{(i-1)/n}\right] \right)
\end{equation*}

Using i), we deduce that 
\begin{equation*}
\mathbb{E}\left[ \left. |b_{1,i}^{n}|^{p}\right\vert \mathcal{F}_{(i-1)/n}%
\right] \leq C\frac{1}{n^{\alpha q+(2\wedge q)/2}}.
\end{equation*}

Now note that%
\begin{equation*}
\{(a_{i}^{n}+u_{i}^{n})a_{i}^{n}<0\}\subset \{|a_{i}^{n}|<|u_{i}^{n}|\}.
\end{equation*}%
Then%
\begin{eqnarray*}
\mathbb{E}\left[ \left. |b_{2,i}^{n}|^{q}\right\vert \mathcal{F}_{(i-1)/n}%
\right] &\leq &C\mathbb{E}\left[ \left. \left(
|a_{i}^{n}|^{q}+|u_{i}^{n}|^{q}\right) \mathbb{I}_{%
\{(a_{i}^{n}+u_{i}^{n})a_{i}^{n}<0\}}\right\vert \mathcal{F}_{(i-1)/n}\right]
\\
&\leq &C\mathbb{E}\left[ \left. |u_{i}^{n}|^{q}\mathbb{I}_{%
\{|a_{i}^{n}|<|u_{i}^{n}|\}}\right\vert \mathcal{F}_{(i-1)/n}\right] .
\end{eqnarray*}%
By H\"{o}lder's inequality,%
\begin{eqnarray*}
\mathbb{E}\left[ \left. |b_{2,i}^{n}|^{q}\right\vert \mathcal{F}_{(i-1)/n}%
\right] &\leq &C\left[ \mathbb{E}\left[ \left. |u_{i}^{n}|^{2q}\right\vert 
\mathcal{F}_{(i-1)/n}\right] \right] ^{1/2}\left[ \Pr \left( \left.
(a_{i}^{n}+u_{i}^{n})a_{i}^{n}<0\right\vert \mathcal{F}_{(i-1)/n}\right) %
\right] ^{1/2} \\
&\leq &C\frac{1}{n^{\alpha q+(2\wedge 2q)/4}}\left[ \Pr \left( \left.
|a_{i}^{n}|<|u_{i}^{n}|\right\vert \mathcal{F}_{(i-1)/n}\right) \right]
^{1/2} \\
&=&C\frac{1}{n^{\alpha q+1/2}}\left[ \Pr \left( \left. \frac{\sqrt{n}}{%
H_{(i-1)/n}}|a_{i}^{n}|<\frac{\sqrt{n}}{H_{(i-1)/n}}|u_{i}^{n}|\right\vert 
\mathcal{F}_{(i-1)/n}\right) \right] ^{1/2}.
\end{eqnarray*}%
Now by Proposition \ref{PropIncMAS}, note that%
\begin{eqnarray*}
&&\Pr \left( \left. \frac{\sqrt{n}}{H_{(i-1)/n}}a_{i}^{n}\leq u\right\vert 
\mathcal{F}_{(i-1)/n}\right) \\
&=&\exp \left( -\frac{1}{\eta _{(i-1)/n}}\left[ e^{-H_{(i-1)/n}n^{-1/2}u}%
\Phi \left( -u+\frac{1}{2}H_{(i-1)/n}n^{-1/2}\right) -\Phi \left( -u-\frac{1%
}{2}H_{(i-1)/n}n^{-1/2}\right) \right] \right) \\
&&\times \Phi \left( u+\frac{1}{2}H_{(i-1)/n}n^{-1/2}\right) .
\end{eqnarray*}%
Let $\lambda _{n}=n^{-1/4}\rightarrow 0$. By Markov's inequality%
\begin{equation*}
\Pr \left( \left. \frac{\sqrt{n}}{H_{(i-1)/n}}|u_{i}^{n}|>\lambda
_{n}\right\vert \mathcal{F}_{(i-1)/n}\right) \leq \frac{n^{1/2}}{H_{(i-1)/n}}%
\frac{\mathbb{E}\left[ \left. |u_{i}^{n}|\right\vert \mathcal{F}_{(i-1)/n}%
\right] }{\lambda _{n}}\underset{n\rightarrow \infty }{\sim }C\frac{1}{%
n^{1/4}}
\end{equation*}%
and moreover%
\begin{equation*}
\Pr \left( \left. \frac{\sqrt{n}}{H_{(i-1)/n}}|a_{i}^{n}|<\lambda
_{n}\right\vert \mathcal{F}_{(i-1)/n}\right) \underset{n\rightarrow \infty }{%
\sim }C\exp \left( -\frac{1}{\eta _{(i-1)/n}}\right) \lambda _{n}\leq C\frac{%
1}{n^{1/4}}\rightarrow 0.
\end{equation*}%
Therefore we have%
\begin{eqnarray*}
&&\Pr \left( \left. \frac{\sqrt{n}}{H_{(i-1)/n}}|a_{i}^{n}|<\frac{\sqrt{n}}{%
H_{(i-1)/n}}|u_{i}^{n}|\right\vert \mathcal{F}_{(i-1)/n}\right) \\
&\leq &\Pr \left( \left. \frac{\sqrt{n}}{H_{(i-1)/n}}|a_{i}^{n}|<\lambda
_{n}\right\vert \mathcal{F}_{(i-1)/n}\right) +\Pr \left( \left. \frac{\sqrt{n%
}}{H_{(i-1)/n}}|u_{i}^{n}|>\lambda _{n}\right\vert \mathcal{F}%
_{(i-1)/n}\right) \\
&\leq &C\frac{1}{n^{1/4}}
\end{eqnarray*}%
for large $n$. Finally, we have 
\begin{equation*}
\mathbb{E}\left[ \left. |b_{2,i}^{n}|^{q}\right\vert \mathcal{F}_{(i-1)/n}%
\right] \leq C\frac{1}{n^{\alpha q+5/8}}
\end{equation*}%
and we deduce that for $q\geq 1$%
\begin{equation*}
\mathbb{E}\left[ \left. |b_{i}^{n}|^{q}\right\vert \mathcal{F}_{(i-1)/n}%
\right] \leq C\left[ \frac{1}{n^{\alpha q+(2\wedge q)/2}}+\frac{1}{n^{\alpha
q+5/8}}\right] \leq C\frac{1}{n^{\alpha q+(q\wedge 5/4)/2}}.
\end{equation*}

b) Study of $w_{i}^{n}$: 
\begin{equation*}
w_{i}^{n}=\sum_{k=1}^{p}C_{p}^{k}(b_{i}^{n})^{k}|a_{i}^{n}|^{p-k}.
\end{equation*}%
For some $q>1$%
\begin{eqnarray*}
&&\mathbb{E}\left[ \left. n^{(p-1)/2}|w_{i}^{n}|\right\vert \mathcal{F}%
_{(i-1)/n}\right] \\
&\leq &\sum_{k=1}^{p}n^{(k-1)/2}C_{p}^{k}\mathbb{E}\left[ \left.
|b_{i}^{n}|^{k}|n^{1/2}a_{i}^{n}|^{p-k}\right\vert \mathcal{F}_{(i-1)/n}%
\right] \\
&\leq &\sum_{k=1}^{p}n^{(k-1)/2}C_{p}^{k}\mathbb{E}\left[ \left.
|b_{i}^{n}|^{qk}\right\vert \mathcal{F}_{(i-1)/n}\right] ^{1/q}\mathbb{E}%
\left[ \left. |n^{1/2}a_{i}^{n}|^{2(p-k)/(1-1/q)}\right\vert \mathcal{F}%
_{(i-1)/n}\right] ^{1-1/q} \\
&\leq &C\sum_{k=1}^{p}n^{(k-1)/2}\frac{1}{n^{\alpha k+(kq\wedge 5/4)/2q}}%
=C\sum_{k=1}^{p}\frac{1}{n^{(\alpha -1/2)k+1/2+(kq\wedge 5/4)/2q}}.
\end{eqnarray*}%
If $1<q<5/4$, we deduce that%
\begin{equation*}
\mathbb{E}\left[ \left. n^{(p-1)/2}|w_{i}^{n}|\right\vert \mathcal{F}%
_{(i-1)/n}\right] =o\left( n^{-1}\right) .
\end{equation*}%
It follows that $\sum_{i=1}^{\left\lfloor nt\right\rfloor -1}\mathbb{E}\left[
\left. n^{(p-1)/2}|w_{i}^{n}|\right\vert \mathcal{F}_{(i-1)/n}\right] 
\overset{u.c.p}{\Longrightarrow }0$ and 
\begin{equation*}
\sum_{i=1}^{\left\lfloor nt\right\rfloor -1}n^{(p-1)/2}w_{i}^{n}\overset{%
u.c.p}{\Longrightarrow }0\text{.}
\end{equation*}

- Step 4) We have%
\begin{equation*}
\sum_{i=1}^{\left\lfloor nt\right\rfloor -1}\zeta
_{i}^{n}=n^{-1/2}\sum_{i=1}^{\left\lfloor nt\right\rfloor -1}\left(
|n^{1/2}\Delta _{i}^{n}\log \eta ^{(i)}|^{p}-m_{p}H_{(i-1)/n}^{p}\right) .
\end{equation*}%
Let%
\begin{equation*}
\kappa _{i}^{n}=\frac{1}{\sqrt{n}}\Psi _{\left\vert \cdot \right\vert
^{p},H_{(i-1)/n}}(\sqrt{n}\Delta _{i}^{n}W_{j},\sqrt{n}\Delta _{i}^{n}W_{k},%
\sqrt{n}Z_{k\backslash j,(i-1)/n})\mathbb{I}_{\{\vee _{l\neq
j,k}Z_{l,(i-1)/n}\leq \wedge _{l=j,k}Z_{l,(i-1)/n}\}}.
\end{equation*}%
By using similar arguments as in Jacod (1998), it is possible to prove that%
\begin{eqnarray*}
\sum_{i=1}^{\left\lfloor nt\right\rfloor }\mathbb{E}\left[ \left. \kappa
_{i}^{n}\right\vert \mathcal{F}_{(i-1)/n}\right] &=&\frac{1}{\sqrt{n}}%
\sum_{i=1}^{\left\lfloor nt\right\rfloor }\mathbb{I}_{\{\vee _{l\neq
j,k}Z_{l,(i-1)/n}\leq \wedge _{l=j,k}Z_{l,(i-1)/n}\}}\varphi
_{p,H_{(i-1)/n}}\left( \sqrt{n}Z_{k/j,(i-1)/n}\right) \\
&&\overset{u.c.p}{\Longrightarrow }\frac{1}{2}\int_{0}^{t}\frac{\lambda
(\varphi _{p,H_{s}})}{H_{s}^{2}}\mathbb{I}_{\{\wedge _{l=j,k}Z_{l,s}>\vee
_{l\neq j,k}Z_{l,s}\}}dL_{Z_{k\backslash j},s}^{0}.
\end{eqnarray*}%
We do not give details, but this can be seen from the following intuitive
arguments%
\begin{eqnarray*}
&&\sum_{i=1}^{\left\lfloor nt\right\rfloor }\mathbb{E}\left[ \left. \kappa
_{i}^{n}\right\vert \mathcal{F}_{(i-1)/n}\right] \\
&\approx &\sqrt{n}\int_{0}^{t}\mathbb{I}_{\{\wedge _{l=j,k}Z_{l,s}>\vee
_{l\neq j,k}Z_{l,s}\}}\varphi _{p,H_{s}}\left( \sqrt{n}Z_{k/j,s}\right) ds \\
&=&\sqrt{n}\int_{0}^{t}\mathbb{I}_{\{\wedge _{l=j,k}Z_{l,s}>\vee _{l\neq
j,k}Z_{l,s}\}}\varphi _{p,H_{s}}\left( \sqrt{n}Z_{k/j,s}\right) \frac{1}{%
2H_{s}^{2}}d\left\langle Z_{k/j,s}\right\rangle \\
&=&\frac{1}{2}\sqrt{n}\int_{\mathbb{R}}\int_{0}^{t}\mathbb{I}_{\{\wedge
_{l=j,k}Z_{l,s}>\vee _{l\neq j,k}Z_{l,s}\}}\varphi _{p,H_{s}}\left( \sqrt{n}%
x\right) \frac{1}{H_{s}^{2}}dL_{Z_{k/j},s}^{x}dx \\
&=&\frac{1}{2}\int_{\mathbb{R}}\int_{0}^{t}\mathbb{I}_{\{\wedge
_{l=j,k}Z_{l,s}>\vee _{l\neq j,k}Z_{l,s}\}}\varphi _{p,H_{s}}\left( z\right) 
\frac{1}{H_{s}^{2}}dL_{Z_{k/j},s}^{z/\sqrt{n}}dz \\
&&\overset{u.c.p}{\Longrightarrow }\frac{1}{2}\int_{\mathbb{R}}\int_{0}^{t}%
\mathbb{I}_{\{\wedge _{l=j,k}Z_{l,s}>\vee _{l\neq j,k}Z_{l,s}\}}\varphi
_{p,H_{s}}\left( z\right) \frac{1}{H_{s}^{2}}dL_{Z_{k/j},s}^{0}dz \\
&=&\frac{1}{2}\int_{0}^{t}\frac{\lambda (\varphi _{p,H_{s}})}{H_{s}^{2}}%
\mathbb{I}_{\{\wedge _{l=j,k}Z_{l,s}>\vee _{l\neq
j,k}Z_{l,s}\}}dL_{Z_{k\backslash j},s}^{0}
\end{eqnarray*}%
where we used the occupation time formula for continuous semimartingales
(see e.g. Proposition 2.1 p 522 in Revuz and Yor (1999)). Moreover $\sqrt{n}%
\sum_{i=1}^{\left\lfloor nt\right\rfloor }\mathbb{E}\left[ \left. (\kappa
_{i}^{n})^{2}\right\vert \mathcal{F}_{(i-1)/n}\right] $ converges u.c.p. to
a non degenerate process. Therefore%
\begin{equation*}
\sum_{i=1}^{\left\lfloor nt\right\rfloor }\kappa _{i}^{n}\overset{u.c.p}{%
\Longrightarrow }\frac{1}{2}\int_{0}^{t}\frac{\lambda (\varphi _{p,H_{s}})}{%
H_{s}^{2}}\mathbb{I}_{\{\wedge _{l=j,k}Z_{l,s}>\vee _{l\neq
j,k}Z_{l,s}\}}dL_{Z_{k\backslash j},s}^{0}.
\end{equation*}

Step 5) By using the same arguments as in the proof of Proposition \ref%
{TCLSimpMS}, we get%
\begin{equation*}
\sum_{i=1}^{\left\lfloor nt\right\rfloor -1}\zeta _{i}^{n}\overset{\mathcal{L%
}-s}{\implies }\tilde{X}_{3,t}
\end{equation*}%
where $\tilde{X}_{3}$ is a process defined on an extension $(\Omega ,%
\mathcal{\tilde{F}},(\mathcal{\tilde{F})}_{t\geq 0},\mathbb{P})$ of $(\Omega
,\mathcal{F},\left( \mathcal{F}\right) _{t\geq 0},\mathbb{P})$, which
conditionally on $\mathcal{F}$ is a continuous Gaussian process, with
independent increments, and whose mean and variance are given respectively by%
\begin{eqnarray*}
\mathbb{\tilde{E}}[\tilde{X}_{3,t}|\mathcal{F]} &=&\frac{1}{2}\sum_{j\geq
1}\sum_{k>j}\int_{0}^{t}\frac{\lambda (\varphi _{p,H_{s}})}{H_{s}^{2}}%
\mathbb{I}_{\{\wedge _{l=j,k}Z_{l,s}>\vee _{l\neq
j,k}Z_{l,s}\}}dL_{Z_{k\backslash j},s}^{0} \\
\mathbb{\tilde{V}}[\tilde{X}_{3,t}|\mathcal{F]} &=&\left(
m_{2p}-m_{p}^{2}\right) \int_{0}^{t}H_{s}^{2p}ds.
\end{eqnarray*}

\end{document}